\documentclass [11pt]{article}
\usepackage[dvips]{graphics}
\usepackage{amsmath}
\usepackage{amssymb}
\usepackage{amscd}
\usepackage{amsthm}
\usepackage{amsopn}
\usepackage{xspace}
\usepackage{verbatim}
\usepackage{amsmath}
\usepackage{amsfonts}
\usepackage{colortbl}
\usepackage{srcltx}
\newtheorem{theorem}{Theorem}
\newtheorem{proposition}{Proposition}[section]
\newtheorem{corollary}{Corollary}[section]
\newtheorem{remark}{Remark}[section]
\newtheorem{lemma}{Lemma}[section]
\numberwithin{equation}{section}

\newcommand\AMSname{AMS subject classifications}
\begin{document}

\title{A Ginzburg-Landau type energy with  weight and with convex potential near zero\\}
\author{Rejeb Hadiji\footnote{Universit\'e Paris-Est, LAMA,
 Laboratoire d'Analyse et de Math\'ematiques Appliqu\'ees, UMR 8050, UPEC, F-94010, Cr\'eteil, France. e-mail: rejeb.hadiji@u-pec.fr  }$\,$ and Carmen
Perugia\footnote{Universit\'a del Sannio, Dipartimento di Scienze e Tecnologie, Via Dei Mulini 59/A Palazzo Inarcassa,
82100, Benevento, Italia. e-mail : cperugia@unisannio.it}}
\date{}

\maketitle

\begin{abstract}
In this paper, we study the asymptotic behaviour of minimizing solutions of a  Ginzburg-Landau type functional with a positive weight and with convex potential near $0$ and we estimate the energy in this case. We also generalize a lower bound for the energy of unit vector field given initially by Brezis-Merle-Rivi\`ere.
\medskip

\medskip

\par
\noindent2010 \AMSname: 35Q56, 35J50,  35B25.
\end{abstract}

\section*{Introduction}
Let $G$ be a bounded, simply connected and smooth domain of $\mathbf{R}^{2}$, 
$g:\partial G\rightarrow S^{1}$ a smooth boundary data of degree $d$ and
$p$ a smooth positive function on $\overline{G}$. We set 
\begin{equation}\label{minp}
p_0 = \min\left\lbrace p(x):x\in \overline {G}\right\rbrace
\end{equation} 
and $\Lambda=p^{-1} (p_0)$.
Let us consider a
$C^{2}$ functional $J:\textbf{R}\rightarrow[0,\infty)$ satisfying
the following conditions :\par
(H1) $J(0)=0$ and $J(t)>0$ on $(0,\infty)$,\par
(H2) $J'(t)>0$ on $(0,1]$,\par
(H3) there exists $\rho_{0}>0$ such that $J''(t)>0$ on $(0,\rho_{0})$.\\
For each $\varepsilon>0$ let $u_{\varepsilon}$ be a minimizer for the following  Ginzburg -Landau type functional
\begin{equation}\label{1.1}
E_{\varepsilon }\left( u\right) =\int_{G}p
\left\vert \nabla u\right\vert ^{2}+\frac{1}{\varepsilon ^{2}}
\int_{G}J\left( 1-\left\vert u\right\vert ^{2}\right)
\end{equation}
defined on the set
\begin{equation}
H_{g}^{1}(G,\mathbf{C})=\left\{ u\in H^{1}(G,\mathbf{C}):\text{ }u=g\text{
on }\partial G\right\}.
\end{equation}
It is easy to prove that $\underset{u\in
H_{g}^{1}(G,\textbf{C})}{\min }E_{\varepsilon}\left(u\right)$ is
achieved by some smooth $u_{\varepsilon}$ which satisfies
\begin{equation}\label{p1}
\left\{
  \begin{array}{cc}
    -div(p\nabla u_{\varepsilon})=\dfrac{1}{\varepsilon^{2}}j(1-|u_{\varepsilon}|^{2})u_{\varepsilon} & $in$\, G \\
    u_{\varepsilon}=g & $on$\,  \partial G, \\
  \end{array}
  \right.
\end{equation}
where $j\left(t\right)=J'\left(t\right)$. In this paper, we are interested in studying the asymptotic behaviour of $u_{\varepsilon}$ and estimate the energy $E_{\varepsilon}(u_{\varepsilon})$ as $\varepsilon\rightarrow 0$ under the assumptions that $p$ has a finite number of minima all lying in $G$ and that it behaves in a "good" way in a neighborhood of each of its minima. More precisely, throughout this paper we shall assume
\begin{equation}\label{P*}
\Lambda=\left\lbrace b_{1},..,b_{N}\right\rbrace\subset G
\end{equation}
and there exist real numbers $\alpha_{k}, \beta_{k}, s_{k}$ satisfying
$0<\alpha_{k} \leq \beta_{k}$ and $s_{k} >  1$ such that
\begin{equation}\label{p}
    \alpha_{k}\left|x-b_{k}\right|^{s_{k}}\leq p(x)-p_{0}\leq\beta_{k}\left|x-b_{k}\right|^{s_{k}}
\end{equation}
in a neighborhood of $b_{k}$ for every $1\leq k\leq N$.\\
The presence of a nonconstant weight function is motivated by the problem of pinning the vortices of $u_\varepsilon$ to some restricted sites, see \cite{DG} and \cite{R} for more detailed physical motivations.
Without loss of generality we assume $d\geq0$. By the way we treat only the case
$d>0$, being the case $d=0$ trivial.\\ The case when  $J(\vert u\vert)=\dfrac{\left(1-|u|^{2}\right)^{2}}{4}$ and $p= {1 \over 2}$ corresponding to the Ginzburg-Landau energy, was studied by several authors since the groundbreaking works of  B\'ethuel-Brezis and H\'elein. More precisely they delt with the case with boundary data satisfying $d=0$ and $d\neq 0$ respectively in \cite{BBH1} and \cite{BBH2}. In this latter work  the case of $G$ starshaped was treated. Eventually in \cite{St}, Struwe gave an argument which works for an arbitrary domain and later del Pino and Felmer in \cite{PF} gave a very simple argument for reducing the general case to the starshaped one. More in particular the method of Struwe is found to be very useful for the case of nonconstant $p$.


In the current paper we will suppose that $card\, \Lambda=N <d$ as this is the more interesting case. Indeed, as already observed in \cite{AS}, singularities of degree $>1$ must occur and in some cases they could be on the boundary. Following the same argument as in \cite{BBH1} or in \cite{AS}, we prove that $u_{\varepsilon_{n}}$ has its zeroes located in $d$ discs, called "bad discs",  with radius $\lambda  \varepsilon_{n}$  where $\lambda > 0$. Outside this discs 
$|u_{\varepsilon_{n}}|$ is close to $1$. For $n$ large each bad discs contains exactly one zero. Thus there are exactly $d_{k}$ zeroes approaching each $b_{k}$ (as $n\rightarrow \infty$). 
In the case $d_{k}>1$ (this must be the case of at least one $k$ if $N<d$), one expects to observe an "interaction energy" between zeroes approaching the same limit $b_{k}$. A complete understanding of this process requires a study of the mutual distances between zeroes of $u_{\varepsilon_{n}}$ which approach the same $b_{k}$. It turns out that these distances depend in a crucial way on the behaviour of the weight function $p$ around its minima points. 
In \cite{AS}, if $s_{k}=2$, it is showed that each $b_{k}$ with $d_{k}>1$ contributes an additional term to the energy, namely $\pi p_{0}\left(d_{k}^{2}-d_{k}\right)\log \left(|\log \varepsilon|^{\frac{1}{2}}\right)$ which is precisely the mentioned interaction energy. In our paper the energy cost of each vortex of degree $>1$ is much less than the previous one.\\

The method of \cite{BBH1,BBH2,St} can be adapted without any difficulties to the case of $J$ satisfying $(H1)-(H3)$ with a zero of finite order at $t=0$. This applies for example to $J(t)=|t|^{k},\,\forall k \geq 2$.  

When $J(\vert u \vert)=\dfrac{\left(1-|u|^{2}\right)^{2}}{4}$ and $p$ not a constant function was studied in \cite{AS,BH1,BH2,BH3}. 
More precisely in \cite{BH1,BH2, BH3} the authors considered the cases $card\,\Lambda=1$ and  $d\geq 1$,  $card\,\Lambda\geq d$ and the case where $p$ has minima on the boundary of the domain. In the first case they highlight  a singularity of degree greater than $1$ when $d> 1$.  More precisely,  if $\Lambda = \{b\} \subset G$, they proved
\begin{equation*}\label{conv2}
u_{\varepsilon_{n}}\rightarrow u_{*}=e^{i\phi} \left(\dfrac{z-b }{|z-b|} \right)^{2d} \quad \hbox {in} \quad C^{1,\alpha}_{loc}\left(\overline{G}\setminus \left\lbrace b  \right\rbrace\right),
\end{equation*}
where $\phi$ is determined by the boundary data $g$.\\

 In the second case, they showed that actually $N=d$, the degree around each $b_{k}$ is equal to $1$ and for a subsequence $\varepsilon_{n}\rightarrow 0$
\begin{equation*}\label{conv1}
u_{\varepsilon_{n}}\rightarrow u_{*}=e^{i\phi}\prod_{j=1}^{d}\dfrac{z-b_{j}}{|z-b_{j}|} \quad \hbox{in}\quad C^{1,\alpha}_{loc}\left(\overline{G}\setminus \left\lbrace b_{1},...,b_{d}\right\rbrace\right),
\end{equation*}
the configuration $\left\lbrace b_{1},...,b_{d}\right\rbrace$ being minimizing for a certain renormalized energy defined in $\Lambda^{d}$. Moreover  they proved the asymptotics $E_{\varepsilon}(u_{\varepsilon})=\pi p_{0}d|\log \varepsilon| + O(1)$.
In the third case, the authors considered the situation when the weight has both minima in the domain and on the boudary.
In \cite{AS}, the authors studied the case $card\, \Lambda<d$ and established the convergence of a subsequence $u_{\varepsilon_{n}}\rightarrow u_{*}$ in $C^{1,\alpha}_{loc}\left(\overline{G}\setminus\left\lbrace b_{1},...,b_{N}\right\rbrace\right)$ for every $\alpha<1$, where the $N$ distinct points $\left\lbrace b_{1},...,b_{N}\right\rbrace$ lie in $\Lambda$ and $u_{*}\in C^{\infty}\left(\overline{G}\setminus\left\lbrace b_{1},...,b_{N}\right\rbrace, S^{1}\right)$ is a solution of 
$$
-div\left(p\nabla u_{*}\right)=p\left|\nabla u_{*}\right|^{2}u_{*}\,\,\,\,in\,\overline{G}\setminus\left\lbrace b_{1},...,b_{N}\right\rbrace,\,\,\,u_{*}=g\,\,on\,\partial G.
$$
Moreover, the degree $d_{k}$ of $u_{*}$ around each $b_{k}$ satisfies $d_{k}\geq 1$ and $\sum_{k=1}^{N}d_{k}=d$.\\

We note that in  \cite{HP} we study the effect of the presence of $\left|u\right|$ in the weight $p(x, u) =p_{0}+s|x|^{k}|u|^{l}$ where $s$ is a small, $k\geq 0$ and $l\geq 0$.

In this article, we are interested in  different types of generalization, starting from the case where  the potential $J$ satisfies $(H1)-(H2)-(H3)$, and $p$ is  non constant. Significative examples are
\begin{equation}
\label{eq:Jk}
J(t)= J_h(t) =\begin{cases}\exp(-1/t^h)& \text{ for }t>0\,,\\
           0 & \text{ for }t\leq 0\,,
\end{cases}
\end{equation}
 for $h>0$. 
 In the present paper, a main new feature is that certain potentials with sufficiently slow growth allow for a vortex energy that is not $\pi|\log\varepsilon| + O(1)$ but instead
 \begin{equation*}
\begin{array}{c}
2\pi p_{0}d_{k}|\log\varepsilon|+2\pi
    p_{0}\frac{d_{k}^{2}-d_{k}}{s_{k}}\log|\log\varepsilon|
-2\pi p_{0}d_{k}I\left(\frac{1}{\varepsilon}\left(|\log\varepsilon|\right)^{-\frac{1}{s_{k}}}\right)+o\left(I\left(\left(|\log\varepsilon| \right)^{\frac{1}{s_{k}}}\right)\right).
\end{array}
\end{equation*}
(see also \cite{HS1,HS2}).
More precisely we want to prove the following result

\begin{theorem}\label{teo1}
For each $\varepsilon>0$, let $u_{\varepsilon }$ be a minimizer for the energy 
\eqref{1.1} over $H_{g}^{1}(G,\mathbf{C})$, with $G$, $g$ as above, $d>0$  and $J$ satisfying (H1)-(H2)-(H3).\\
i) For a subsequence $\varepsilon_{n}\rightarrow 0$ we have
\begin{equation}\label{conv1}
u_{\varepsilon_{n}}\rightarrow u_{*}=e^{i\phi}\prod_{j=1}^{N} ( \dfrac{z-b_{j}}{|z-b_{j}|} )^{d_j}\,\,\,\,in\, C^{1,\alpha}_{loc}\left(\overline{G}\setminus \left\lbrace b_{1},...,b_{N}\right\rbrace\right)
\end{equation} 
for every $\alpha<1$, where the $N$ distinct points $\left\lbrace b_{1},...,b_{N}\right\rbrace$ 
lie in $\Lambda$, $\Sigma_{j = 1}^{N} d_j = d$  and $\phi$ is a smooth harmonic function determined by the requirement $u_{*}=g$ on $\partial G$.\\
ii) Setting
 \begin{equation*}\label{HS2.11}
     I(R) = \dfrac{1}{2}\displaystyle\int_{\frac{1}{R^{2}}}^{j\left(\rho_{0}\right)}\frac{j^{-1}(t)}{t}dt\,
\end{equation*}
we have
\begin{equation}\label{conv2}
\begin{split}
                               E_{\varepsilon_{n}}\left(u_{\varepsilon_{n}}\right)    =&  2\pi p_{0}d\log\frac{1}{\varepsilon_{n}}+2\pi
    p_{0}\left(\displaystyle\Sigma_{k=1}^{N}\frac{d_{k}^{2}-d_{k}}{s_{k}}\right)\log\log\frac{1}{\varepsilon_{n}}
\\
&-2\pi
p_{0}dI\left(\dfrac{1}{\varepsilon_{n}}\left(\log\frac{1}{\varepsilon_{n}}\right)^{-\frac{1}{s_{k}}}\right)+o\left(I\left(\left(\log\frac{1}{\varepsilon_{n}}\right)^{\frac{1}{s_{k}}}\right)\right).
\end{split}
\end{equation}
\end{theorem}
As it is showed in \cite{HS1}, $ \lim_{R\rightarrow \infty}\dfrac{I(R)}{\log R}=0$ hence the leading term in the energy is always of order $o(|\log \varepsilon|)$. Moreover it is easy to see that $I(R)$ is a positive, monotone increasing, concave function of $\log R$ for $R$ large (see \cite{HS1}).
The proof of Theorem \ref{teo1} consists of two main ingredients: the method of Struwe \cite{St}, as used also in \cite{AS} in order to locate the "bad discs", (i. e. a finite collection of discs of radius $O(\varepsilon)$ which cover the set $\left\lbrace x: |u_{\varepsilon}(x)<\frac{1}{2}|\right\rbrace$) and the generalization of a result of Brezis, Merle and Rivi\`{e}re \cite{BMR} which will play an important role in finding the lower bound of the energy. More precisely in Theorem \ref{teoBMR.4}, we will bound from below the energy of a regular map defined away from 
some points  $a_{1},a_{2},...,a_{m}$  in $B_{R}(0)$ such that
 $0<a\leq\left|u\right|\leq 1$ in $\Omega$, $deg\left(u,\partial B_{R}(a_{j}\right)=d_{j}$
and with a bound potential  by using the reference map $u_{0}(z)=\left(\frac{z-a_{1}}{\left|z-a_{1}\right|}\right)^{d_{1}}\left(\frac{z-a_{2}}{\left|z-a_{2}\right|}\right)^{d_{2}}....\left(\frac{z-a_{m}} {\left|z-a_{m}\right|}\right)^{d_{m}}$. 
After the results of \cite{BMR},  Han and  Shafrir , Jerrard, Sandier, Struwe obtained the essential lower bounds for the Dirichlet energy of a unit vector field, see \cite{HaS}, \cite{J}, \cite{S} and \cite{St}.
The paper is organized as follows. In Section 1, we recall some definitions and results contained in \cite{HS1}. Section 2 is devoted to prove the generalization of Theorem 4 of \cite{BMR} which will be useful for obtaining a precise lower bound of the energy for our case. In Section 3 we prove Theorem \ref{teo1} by stating an upper and a lower bound for the energy $(\ref{1.1})$. Finally, as a corollary of upper and lower bounds of the energy, we find an estimate of the mutual distances between bad discs approaching the same singularity  $b_k$.

\section{Recalls}
In this section we recall some results proved in \cite{HS1} useful in the sequel. Let us consider the following quantity, introduced in \cite{HS1} which will play an important role in our study
\begin{equation}\label{1.2}
    I\left(R,c\right)=\sup\left\{\int_{1}^{R}\frac{1-f^{2}}{r}dr:\int_{1}^{R}J\left(1-f^{2}\right)rdr\leq c\right\}
\end{equation}
for any $R>1$ and $c>0$.\\
\begin{lemma}\label{2.1HS}
For every $R>0$ and $c>0$, there exists a maximizer $f_{0}=f_{0}^{(R)}$ in \eqref{1.2} satisfying $0\leq f_{0}(r)\leq 1$ for every $r$ such that $f_{0}(r)$ is nondecreasing. Moreover, if $r_{0}=r_{0}(c)$ is defined by the equation
$$
c=J(1)\left(\dfrac{r_{0}^{2}-1}{2}\right),
$$
then there exists $\widetilde{r_{0}}=\widetilde{r_{0}}(c,R)\in [1,r_{0}]$ such that
$$
f_{0}(r)\left\lbrace
\begin{array}{ll}
=0\,\,\,\,\,\,\,if\,r\in[1,R]\,and\,r<\widetilde{r_{0}},\\
>0 \,\,\,\,\,\,\,if\,\,r>\widetilde{r_{0}}.
\end{array}
\right.
$$
Furthermore
$$
\int_{1}^{R}J\left(1-f_{0}^{2}\right)rdr=c,\,\,\,\forall R>r_{0}
$$
and
$$
j\left(1-f_{0}^{2}\right)=\dfrac{1}{\lambda r^{2}},\,\,r>\widetilde{r_{0}}
$$
for some $\lambda=\lambda(R,c)>0$.
\end{lemma}
Moreover it holds
\begin{lemma}\label{2.2HS}
There exist two constants $\kappa_1 >0, \kappa_2 >0$ such that
\begin{equation}\label{errata1HS}
\kappa_1 min(1, \frac{1}{c}) \leq \lambda \leq \kappa_2 (1 + \frac{1}{c}),\,\,\,R\geq r_{0}+1.
\end{equation}
\end{lemma}
\noindent 
Actually, the proof of the previous lemma shows that the estimate of $\lambda$ is uniform for $c$ lying in a bounded interval.

\begin{lemma}\label{2.3HS}
For every $c > 1$ there exists a constant $C( c)$ 
such that for every $c _{1},c_{2} \in [1 /c , c]$ we have
\begin{equation}\label{errata2HS}
\left|I(R,c_{1})-I(R,c_{2})\right|\leq C(c) \,\,\,\, \forall R\geq 1.
\end{equation}
\end{lemma}
\noindent In view of Lemma \ref{2.3HS} it is natural to set
$$
I(R)=I(R,1)
$$
and for any fixed $c_0 > 1$ we have  
\begin{equation}\label{errata3HS}
\left|I(R,c)-I(R)\right|\leq   C(c_0),\,\,\,\forall c\in [ 1/c_0, c_0] \,\,\forall R\geq 1.
\end{equation}
We  recall some properties of $I(R)$.
\begin{lemma}\label{2.4HS}
We have 
\begin{equation}\label{HS2.11}
     I(R) = \dfrac{1}{2}\displaystyle\int_{\frac{1}{R^{2}}}^{j\left(\eta_{0}\right)}\frac{j^{-1}(t)}{t}dt\,\,\,\,
    \forall R\geq 1.
\end{equation}
In particular,
\begin{equation}\label{HS2.12}
    \lim_{R\rightarrow \infty}\frac{I(R)}{\log R}=0.
\end{equation}
\end{lemma}
Moreover for every  $\alpha>0$ there
exists a constant $C_{1}\left(\alpha\right)$ such that
\begin{equation}\label{HS2.13}
    \left|I\left(\alpha R\right)-I(R)\right|\leq C_{1}\left(\alpha\right)
\end{equation}
for $R>\max\left(1,\dfrac{1}{\alpha}\right)$ and
$c\in\left(0,c_{0}\right]$.\\
The next lemma provides an estimate we shall use in the proof of the upper bound in subsection \ref{sub2.1}.
\begin{lemma}\label{2.5HS}
We have
$$
\int_{\mu_{0}}^{R}\left(f_{0}'\right)^{2}\leq C,\,\,\forall R>\mu_{0}
$$
where $\mu_{0}=\max\left(r_{0}(1), \dfrac{1}{\sqrt{a j\left(\rho_{0}\right)}}\right)$ being 
$r_{0}(1)$ and $a$ defined respectively as in Lemma \ref{2.1HS} and Lemma \ref{2.2HS}.
\end{lemma}
In Theorem 1 we will need a similar functional to that of (\ref{1.2}). Hence for $R>1$ and $c> 0$ we set
\begin{equation}\label{1.2*}
  \widetilde{I} \left(R,c\right)=\sup\left\{\int_{1}^{R}\frac{1-f^{2}}{r}dr+4\int_{1}^{R}\frac{\left( 1-f^{2}\right) ^{2}}{r}dr:\int_{1}^{R}J\left(1-f^{2}\right)rdr\leq c\right\}.
\end{equation}
\\
Now, let us  recall an important relation between the two functionals (\ref{1.2}) and (\ref{1.2*}).
\begin{lemma}\label{2.6HS}
There exists a constant $C=C(c)$ such that
\begin{equation}\label{HS2.29}
    \left|\widetilde{I}\left(R,c\right)-I(R,c)\right|\leq C
\end{equation}
for $R  >1$.
\end{lemma}

\begin{lemma}\label{2.7HS}
There exists a constant $\kappa$ such that for every  $c > 0$, $\alpha$ $>0$, 
$$
\left|I\left(\alpha R,c\right)-I(R)\right|\leq \kappa (c_{0},\alpha)
$$
$$
\left|\widetilde{I}\left(\alpha R,c\right)-I(R)\right|\leq C_{1}(c_{0},\alpha)
$$
for $R> max \left(1, \dfrac{1}{\alpha}\right)$ and $c\in (0,c_{0}]$.
\end{lemma}

The next two propositions, dealing with a lower bound for the energy in a simple annulus and in a more general perforated domain respectively, will play an important role in the proof of our lower bound stated in subsection \ref{sub3.2} (see \cite{HS1} for ditails).
\begin{proposition}\label{3.2HS}
Let $A_{R_{1},R_{2}}$ denotes the annulus $\left\lbrace R_{1}< |x|< R_{2}\right\rbrace$ and let $$u\in C^{1}\left(A_{R_{1},R_{2}},\mathbb{C}\right)\cap C\left(\overline{{A_{R_{1},R_{2}}}},\mathbb{C}\right)$$ satisfy
$$
\deg\left(u, \partial B_{R_{j}}(0)\right)=d,\,\,\,j=1,2,
$$
\medskip
$$
\dfrac{1}{2}\leq|u|\leq 1\,\,on\,A_{R_{1},R_{2}}
$$
and
$$
\dfrac{1}{R_{1}^{2}}\int_{A_{R_{1},R_{2}}}J\left(1-|u|^{2}\right)\leq c_{0},
$$
for some constant $c_{0}$. Then there exists a constant $c_{1}$ depending only on $c_{0}$ such that
$$
\int_{A_{R_{1},R_{2}}} |\nabla u|^{2} \geq 2\pi d^{2}\left(\log\dfrac{R_{2}}{R_{1}}-I\left(\dfrac{R_{2}}{R_{1}}\right)\right)-d^{2}c_{1}.
$$
\end{proposition}

\begin{proposition}\label{3.3HS}
Let $x_{1},x_{2},...,x_{m}$ be $m$ points in $B_{\sigma}(0)$ satisfying
$$
|x_{i}-x_{j}|\geq 4\delta, \forall i\neq j\,\,\,and\,\,|x_{i}|<\dfrac{\sigma}{4},\,\,\forall i,
$$
with $\delta\leq\dfrac{\sigma}{32}$. Set $\Omega=B_{\sigma}(0)\setminus \bigcup _{j=1}^{m} B_{\delta}(x_{j)}$ and let $u$ be a $C^{1}$-map from $\Omega$ into $\mathbb{C}$, which is continuous on $\partial \Omega$ satisfying
$$
\deg\left(u, \partial B_{\sigma}(x_{j})\right)=d_{j},\,\,\,\forall j
$$
$$
\dfrac{1}{2}\leq|u|\leq 1\,\,in\,\Omega
$$
and
$$
\dfrac{1}{\delta^{2}}\int_{\Omega}J\left(1-|u|^{2}\right)\leq K.
$$
Then, denoting $d=\sum_{j=1}^{m}d_{j}$, we have
$$
\int_{\Omega} |\nabla u|^{2} \geq 2\pi |d|\left(\log\dfrac{\sigma}{\delta}-I\left(\dfrac{\sigma}{\delta}\right)\right)-C
$$
with $C=C\left(K,m,\sum_{j=1}^{m}|d_{j}|\right)$.
\end{proposition}
\section{Lower bound for the energy of unit vector fields }
In this section we will  generalize Theorem 4 of \cite{BMR}. To this aim let $a_{1},a_{2},...,a_{m}$ be $m$ points in $B_{R}(0)$ such that
\begin{equation}\label{BMR1}
\left|a_{i}-a_{j}\right|\geq 4 R_{0},\quad \forall i\neq j
\end{equation}
and 
\begin{equation}\label{BMR1i}
\left|a_{i}\right|\leq \frac{R}{2},\quad \forall i,
\end{equation}
with 
\begin{equation}\label{BMR1ii}
R_{0}\leq\frac{R}{4}.
\end{equation} 
Set 
$$\Omega = B_{R}(0)\backslash\bigcup_{j=1}^{m}B_{R_{0}}(a_{j})$$
and let $u$ be a $C^{1}-map$ from $\Omega$ into $\mathbf{C}$ which is continuous
on $\partial\Omega$.\\
We suppose that  
\begin{equation}\label{BMR2}
    0<a\leq\left |u\right|\leq 1 \,\text{in} \, \Omega
\end{equation}
and
\begin{equation}\label{BMR3}
    \frac{1}{R_{0}^{2}}\int_{\Omega}J\left(1-\left|u\right|^{2}\right)\leq K,
\end{equation}
for some constants $a$ and $K$.\\
Let us observe that \eqref{BMR2} implies
$$
\text{deg}\left(u,\partial B_{R}(a_{j}\right))=d_{j}\,\, \quad\forall j$$
is well defined. Hence, let us denote $d=\sum _{j=1}^{m}\left|d_{j}\right|$ and consider the map
\begin{equation}\label{BMR4}
    u_{0}(z)=\left(\frac{z-a_{1}}{\left|z-a_{1}\right|}\right)^{d_{1}}\left(\frac{z-a_{2}}{\left|z-a_{2}\right|}\right)^{d_{2}}....\left(\frac{z-a_{m}}
    {\left|z-a_{m}\right|}\right)^{d_{m}}.
\end{equation}
We want to prove the following result
\begin{theorem}\label{teoBMR.4}
Let us suppose that \eqref{BMR1}$\div$\eqref{BMR3} hold, then we have
\begin{equation}\label{BMR5}
\begin{array}{c}
\displaystyle\int_{\Omega}p\left|\nabla u\right|^{2}\geq p_{0}\int_{\Omega}\left|\nabla u_{0}\right|^{2}-2\pi p_{0}\left(\sum_{i=1}^{m}d_{i}^{2}\right)I\left(\dfrac{R}{R_{0}}\right)+\\
\\
\displaystyle-2\pi \left(1-a^{2}\right)p_{0}\sum_{i\neq j}\left|d_{i}\right|\left|d_{j}\right|\log \frac{R}{|a_{i}-a_{j}|}-C,
\end{array}\end{equation}
where $C$ is a constant depending only on $p_{0}$, $a$, $d$, $m$ and $K$.
\end{theorem}
\textbf{Proof} Let us set $\rho=\left|u\right|$ so that $u=\rho e^{i\varphi}$ locally in $\Omega$.
Hence we have $$\left|\nabla u\right|^{2}=\left|\nabla \rho\right|^{2}+\rho^{2}\left|\nabla \varphi\right|^{2}.$$
Similarly, we can set $u_{0}=e^{i\varphi_{0}}$ locally in $\Omega$ which implies $\left|\nabla u_{0}\right|=\left|\nabla \varphi_{0}\right|$ and 
\begin{equation}\label{gradfi0}
    \nabla \varphi_{0}(z)=\sum_{i=1}^{m}d_{i}\frac{V_{i}(z)}{\left|z-a_{i}\right|},
\end{equation}
where \begin{equation*}
    V_{i}(z)=\left(-\frac{y-a_{i}}{\left|z-a_{i}\right|};\frac{x-a_{i}}{\left|z-a_{i}\right|}\right)
\end{equation*} 
is the unit vector tangent to the circle of radius $\left|z-a_{i}\right|$ centered at $a_{i}$.

By introducing the function $\psi=\varphi-\varphi_{0}$, we can write $u=\rho u_{0}e^{i\psi}$ and have
\begin{equation}\label{BMR6}
    \left|\nabla u\right|^{2}=\left|\nabla \rho\right|^{2}+\rho^{2}\left|\nabla \varphi_{0}+\nabla\psi\right|^{2}.
\end{equation}
By \eqref{minp} and $\left( \ref{BMR6}\right)$ we get
\begin{equation*}
\int_{\Omega}p\left|\nabla u\right|^{2}\geq p_{0}\int_{\Omega}\left|\nabla\rho\right|^{2}+
    p_{0}\int_{\Omega}\rho^{2}\left|\nabla\varphi_{0}\right|^{2}
+p_{0}\int_{\Omega}\rho^{2}\left|\nabla\psi\right|^{2}
    +2p_{0}\int_{\Omega}\rho^{2}\nabla\varphi_{0}\nabla\psi.
\end{equation*}
By adding and subtracting one in the second and fourth integral and by $\left( \ref{BMR2}\right)$, we get
\begin{equation}\label{BMR33*}
\begin{split}
\int_{\Omega}p\left|\nabla u\right|^{2}\geq &-p_{0}\int_{\Omega}\left(1-\rho^{2}\right)\left|\nabla\varphi_{0}\right|^{2}
                  +p_{0}\int_{\Omega}\left|\nabla\varphi_{0}\right|^{2}+p_{0}a^{2}\int_{\Omega}\left|\nabla\psi\right|^{2}
                  \\
                  &+2p_{0}\int_{\Omega}\left(\rho^{2}-1\right)\nabla\varphi_{0}\nabla\psi+
                  2p_{0}\int_{\Omega}\nabla\varphi_{0}\nabla\psi.
\end{split}
\end{equation}
Using $  2AB\geq -|A|^{2}-|B|^{2}$,
for $A=2\left(\rho^{2}-1\right)\nabla\varphi_{0}$ and $B=\dfrac{\nabla\psi}{2}$, we can write
\begin{equation}\label{BMR7}
    \begin{split}
        \int_{\Omega}p\left|\nabla u\right|^{2}\geq &-p_{0}\int_{\Omega}\left(1-\rho^{2}\right)\left|\nabla\varphi_{0}\right|^{2}+p_{0}\int_{\Omega}\left|\nabla\varphi_{0}\right|^{2}+p_{0}a^{2}\left\|\nabla\psi\right\|_{2}^{2}
     \\
        &  -4p_{0}\int_{\Omega}\left(\rho^{2}-1\right)^{2}\left|\nabla\varphi_{0}\right|^{2}-\dfrac{p_{0}}{4}\left\|\nabla\psi\right\|^{2}_{2}+2p_{0}\int_{\Omega}\nabla\varphi_{0}\nabla\psi.
     \end{split}
\end{equation}
As in Theorem 4 of [8] it holds
\begin{equation}\label{BMR8}
    \left|\int_{\Omega}\nabla\varphi_{0}\nabla\psi\right|\leq Cm|d|\left\|\nabla\psi\right\|_{2},
\end{equation}
for some universal constant $C$, hence $\left( \ref{BMR7}\right) $ becomes
\begin{equation}\label{BMR9}
    \begin{split}
      \int_{\Omega}p\left|\nabla u\right|^{2}\geq &\, p_{0}\int_{\Omega}\left|\nabla\varphi_{0}\right|^{2}-\left[ p_{0}\int_{\Omega}\left(1-\rho^{2}\right)\left|\nabla\varphi_{0}\right|^{2}+4p_{0}\int_{\Omega}\left(\rho^{2}-1\right)^{2}\left|\nabla\varphi_{0}\right|^{2}\right]\\
        &+p_{0}\left(a^{2}-\frac{1}{4}\right) \left\|\nabla\psi\right\|^{2}_{2}-2 p_{0}Cm|d|\left\|\nabla\psi\right\|_{2}.
    \end{split}
\end{equation}
Now let us denote $X=\left\|\nabla\psi\right\|_{2}$ and consider the following function
\begin{equation*}
    Y=\left(a^{2}-\frac{1}{4}\right)X^{2}-2Cm\left|d\right|X.
\end{equation*}
If $a>\dfrac{1}{2}$, it reaches its minimum value $Y_{min}=-\dfrac{C^{2}m^{2}\left|d\right|^{2}}{a^{2}-\frac{1}{4}}$ at $X_{min}=\dfrac{Cm\left|d\right|}{a^{2}-\frac{1}{4}}$. Then we get
\begin{equation}\label{BMR10}
\int_{\Omega}p\left|\nabla u\right|^{2}\geq p_{0}\int_{\Omega}\left|\nabla u_{0}\right|^{2}-p_{0}\left[\int_{\Omega}\left(1-\rho^{2}\right)\left|\nabla\varphi_{0}\right|^{2}+
        4\int_{\Omega}\left(\rho^{2}-1\right)^{2}\left|\nabla\varphi_{0}\right|^{2}\right]-C
\end{equation}
where $C$ is a constant depending only on $p_{0}$, $a$, $d$ and $m$. 

Taking into account  (\ref{1.2}), (\ref{1.2*}) and (\ref{HS2.29}), in order to get our result, it is enough to estimate the following term
\begin{equation}\label{BMRI}
\int_{\Omega}\left(1-\rho^{2}\right)\left|\nabla\varphi_{0}\right|^{2}.
\end{equation}
To this aim let us observe that \eqref{gradfi0} implies
\begin{equation*}
\left|\nabla\varphi_{0}(z)\right|^{2} \leq\sum_{i=1}^{m} \frac{d_{i}^{2}}{\left|z-a_{i}\right|^{2}}+
\sum_{i\neq j}\frac{d_{i}d_{j}}{\left|z-a_{i}\right|\left|z-a_{j}\right|}.
\end{equation*}
Then $\left(\ref{BMRI}\right)$ can be written as
\begin{equation}\label{BMR11}
\begin{split}
\int_{\Omega}\left(1-\rho^{2}\right)\left|\nabla\varphi_{0}\right|^{2}=& \int_{\Omega}\left(1-\rho^{2}\right)\left[\sum_{i=1}^{m}\frac{d_{i}^{2}}{\left|z-a_{i}\right|^{2}}+
\sum_{i\neq j}\frac{d_{i}d_{j}}{\left|z-a_{i}\right|\left|z-a_{j}\right|}\right]\\
\leq & \sum_{i=1}^{m}  d_i^2 \int_{\Omega} {1-\rho^2\over | z-a_{i} |^{2}}+
\sum_{i\neq j}d_{i}d_{j}\int_{\Omega}\frac{1 -\rho^2}{\left|z-a_{i}\right|\left|z-a_{j}\right|}\\
=& \sum_{i=1}^{m}d_{i}^{2}A_{i}+B.
\end{split}
\end{equation}
\noindent
Let us analyze each term separately. In order to estimate $A_i$ for every $i=1,..,m$,  let us introduce $\delta_{i}=\text{dist} \left(a_{i},\partial B_{R}(0)\right)$ and observe that $\dfrac{R}{2}\leq\delta_{i}\leq R$ as a consequence of $(\ref{BMR1ii})$.\\
Therefore for any fixed $i$, by definition $(\ref{1.2*})$, it holds
\begin{equation}\label{a}
A_{i}= \int_{\Omega}\frac{1-\rho^{2}}{\left|z-a_{i}\right|^{2}}dz\leq \int_{ B_{R}(0)\setminus B_{R_{0}}(a_{i})}\frac{1-\rho^{2}}{\left|z-a_{i}\right|^{2}}dz\leq 2\pi I\left(\dfrac{\delta_{i}}{R_{0}}\right)\leq 2\pi I\left(\dfrac{R}{R_{0}}\right)+C
\end{equation}
where $C$ {depends only on $K$ defined in \eqref{BMR3}} but is independent of $R$, $R_{0}$ and $a_{i}$. For the second term, acting as in Theorem 5 of \cite{BMR} and using $(\ref{BMR2})$ we obtain \begin{equation}\label{b}
|B|\leq \sum_{i\neq j}\left|d_{i}\right|\left|d_{j}\right|\int_{\Omega}\frac{1 -\rho^2}{\left|z-a_{i}\right|\left|z-a_{j}\right|}dz \leq 2\pi(1-a^{2}) \sum_{i\neq j}\left|d_{i}\right|\left|d_{j}\right|\log \dfrac{R}{\left|a_{i}-a_{j}\right|}+C.
\end{equation}
where $C$ depends only on $m$ and $d$.\\

Then by putting together $(\ref{a})$ and $(\ref{b})$ into $(\ref{BMR11})$
we get 
\begin{equation}\label{BMR12}
   \int_{\Omega}\left(1-\rho^{2}\right)\left|\nabla\varphi_{0}\right|^{2}\leq 2\pi \left(\sum_{i=1}^{m}  d_{i}^2\right) I\left(\frac{R}{R_{0}}\right)+2\pi(1-a^{2}) \sum_{i\neq j}\left|d_{i}\right|\left|d_{j}\right|\log \dfrac{R}{\left|a_{i}-a_{j}\right|}+C.
\end{equation}
where $C$ depending on $ K$, $a$, $m$ and  $d$ but doesn't depend on $R$, $R_{0}$ and $a_{i}$ for every $i=1,..,m$.

Finally, by $\left( \ref{BMR10}\right) $ and $\left( \ref{BMR12}\right) $ we get $\left( \ref{BMR5}\right) $.\\

Under the same hypotheses of Theorem \ref{teoBMR.4}, as an immediate consequence of $(\ref{BMR5})$ and Theorem 5 of \cite{BMR}, we get the following result

\begin{corollary}\label{cor1.1} Let us suppose that \eqref{BMR1}$\div$\eqref{BMR3} hold, then we have
\begin{equation}\label{BMR5*}
\begin{split}
\displaystyle\int_{\Omega}p\left|\nabla u\right|^{2}\geq &\, 2\pi p_{0} \left(\sum_{i=1}^{m}d_{i}^{2}\right) \left( \log \dfrac{R}{R_{0}} -    I\left(\dfrac{R}{R_{0}} \right)\right)\\
 &+2\pi p_{0}\displaystyle\sum_{i\neq j}\left(-\left(1-a^{2}\right)\left|d_{i}\right|\left|d_{j}\right|+d_{i}d_{j}\right)\log \frac{R}{|a_{i}-a_{j}|}-C,
\end{split}
\end{equation}
where $C$ is a constant depending only on $p_{0}$, $a$, $d$, $m$ and $K$.\\
\end{corollary}

\begin{remark}\label{positif}
If $d_{i}\geq 0$ for $i = 1, ..., m$ then \eqref{BMR5*} becomes

\begin{equation}\label{BMR5*positif}
\displaystyle\int_{\Omega}p\left|\nabla u\right|^{2}\geq 2\pi p_{0} \left(\sum_{i=1}^{m}d_{i}^{2}\right) \left( \log \dfrac{R}{R_{0}} -    I\left(\dfrac{R}{R_{0}} \right)\right)+2\pi p_{0}a^{2}\displaystyle\sum_{i\neq j}d_{i}d_{j}\log \frac{R}{|a_{i}-a_{j}|}-C,
\end{equation}
where $C$ is a constant depending only on $p_{0}$, $a$, $d$, $m$ and $K$.
\end{remark}

\section{Proof of Theorem \ref{teo1}}
Throughout this section, for any subdomain $D$ of $G$ we shall denote
\begin{equation}
E_{\varepsilon}\left(u,D\right)=\int_{D}p|\nabla u|^{2}+\dfrac{1}{\varepsilon^{2}}\int_{D}J\left(1-|u|^{2}\right)
\end{equation}
and if $D=G$ we simply write $E_{\varepsilon}(u)$. Moreover, similarly to Proposition \ref{3.2HS}, we will use the following notation 
\begin{equation}\label{annulus}
B_{R_{1},R_{2}}(b) =\left\lbrace R_{1}< |x-b|< R_{2}\right\rbrace
\end{equation}
for the annulus centered in $b$ and with radius $R_1$ and $R_2$. \\
Our main result of this section is the asymptotic behavior of the energy for minimizers which will give \eqref{conv2} of Theorem \ref{teo1}. More precisely we prove the following result
\begin{proposition}\label{as}
Assume $\eqref{P*}$ and $\eqref{p}$ hold true. Then for a subsequence $\varepsilon_{n}\rightarrow 0$ we have
\begin{equation}\label{up11*}
\begin{split}
   \displaystyle                            E_{\varepsilon_{n}}\left(u_{\varepsilon_{n}}\right)    =&\,  2\pi p_{0}d\log\frac{1}{\varepsilon_{n}}+2\pi
    p_{0}\left(\displaystyle\Sigma_{k=1}^{N}\frac{d_{k}^{2}-d_{k}}{s_{k}}\right)\log\log\frac{1}{\varepsilon_{n}}
\\
&-2\pi
p_{0}dI\left(\dfrac{1}{\varepsilon_{n}}\left(\log\frac{1}{\varepsilon_{n}}\right)^{-\frac{1}{s_{k}}}\right)+o\left(I\left(\left(\log\frac{1}{\varepsilon_{n}}\right)^{\frac{1}{s_{k}}}\right)\right).
\end{split}
\end{equation}
\end{proposition}

\subsection{An upper bound for the energy}\label{sub2.1}

Let us prove an upper bound for the functional $\left(\ref{1.1}\right)$.
\begin{proposition}\label{upb}
Let us suppose that $\eqref{P*}$ and $\eqref{p}$ hold true. Then for a subsequence $\varepsilon_{n}\rightarrow 0$ we have 
\begin{equation}\label{up11***}
\begin{split}
                               E_{\varepsilon_{n}}\left(u_{\varepsilon_{n}}\right)\leq &2\pi p_{0}d\log\frac{1}{\varepsilon_{n}}+2\pi
    p_{0}\left(\displaystyle\Sigma_{k=1}^{N}\frac{d_{k}^{2}-d_{k}}{s_{k}}\right)\log\log\frac{1}{\varepsilon_{n}}\\
&-2\pi
p_{0}d I\left(\frac{1}{\varepsilon_{n}}\left(\log\frac{1}{\varepsilon_{n}}\right)^{-\frac{1}{s_{k}}}\right)+O(1).
\end{split}
\end{equation}
\end{proposition}

\textbf{Poof.}
Let $\eta_{0}>0$ satisfy
\begin{equation*}
    0<\eta_{0}<\frac{1}{4}{\min }\left(\underset{i\neq j}{\min }\left|\bar{b}_{i}-\bar{b}_{j}\right|,\underset{i=1,..,N}{\min }dist\left(\bar{b}_{i},\partial
    G\right)\right)
\end{equation*}
and fix $k=1,...,N$. Set
\begin{equation}\label{definitionTepsilon}
T_{\varepsilon_{n}}=\left(\log\frac{1}{\varepsilon_{n}}\right)^{-\frac{1}{s_{k}}}.
\end{equation}

We  will construct a function
$U_{\varepsilon_{n}}(x)$ defined in $\bigcup_{k=1}^{N}B_{\eta_{0}}\left(\bar{b}_{k}\right)$. From this point onwards the proof will develop into three steps. \\

\textbf{Step 1.}  

We define
\begin{equation}\label{refmap4}
    U_{\varepsilon_{n}}(x)=U_{\varepsilon_{n}}^{k}(x)=
    \left(\dfrac{x-\bar{b}_{k}}{|x-\bar{b}_{k}|}\right)^{d_{k}}
    \,\,\,\,\,\,\,\,\,\,\,\,\,\,\,\,\,\,\,\,\,\text{ on }\, B_{T_{\varepsilon_{n}}, \eta_{0}}(\bar{b}_{k}).  
\end{equation}
By following a similar argument as in \cite{AS}, it is  easy to show that 

\begin{equation}\label{refmap1}
  E_{\varepsilon_{n}}\left(U_{\varepsilon_{n}}^{k},B_{\eta_{0}}(\bar{b}_{k})\setminus
    \overline{B_{T_{\varepsilon_{n}}}(\bar{b}_{k})}\right)\leq 2\pi p_{0}\frac{d_{k}^{2}}{s_{k}}\log\log\frac{1}{\varepsilon_{n}}+O(1).
\end{equation}
\textbf{Step 2.} 
Let us fix $d_{k}$ equidistant points $x_{1}^{n},x_{2}^{n},...,x_{d_{k}}^{n}$ on the circle $\partial B_{\frac{T_{\varepsilon_{n}}}{2}}\left(\bar{b}_{k}\right)$ and set
$$A_{\varepsilon_{n}}=B_{T_{\varepsilon_{n}}}(\bar{b}_{k},)\setminus
\bigcup_{j=1}^{d_{k}}B_{\frac{T_{\varepsilon_{n}}}{10d_{k}}}\left(x_{j}\right).$$
 We define $U_{\varepsilon_{n}}$ 

 as an $S^{1}$-valued map which minimizes the
energy $\displaystyle\int_{A_{\varepsilon_{n}}}p\left|\nabla u\right|^{2}$ among
$S^{1}$-valued maps for the boundary data
$\left(\dfrac{x-\bar{b}_{k}}{\left|x-\bar{b}_{k}\right|}\right)^{d_{k}}$ on
$\partial B_{T_{\varepsilon_{n}}}\left(\bar{b}_{k},\right)$ and
$\dfrac{x-x_{j}}{\left|x-x_{j}\right|}$
on $\partial B_{\frac{T_{\varepsilon_{n}}}{10d_{k}}}\left(x_{j},\right)$,
$j=1,..,d_{k}$. 
Clearly we have
\begin{equation}\label{refmap2}
     E_{\varepsilon_{n}}\left(U_{\varepsilon_{n}}(x),A_{\varepsilon_{n}}\right)\leq C.
\end{equation}
Now, let us fix $j\in \left\lbrace 1,..,d_{k}\right\rbrace$, let $\vartheta_{j}$ denote a polar coordinate around $x_{j}$ and let $f_{0}(r)$ be a maximizer for
$I\left(\frac{1}{\varepsilon_{n}}\left(\log\frac{1}{\varepsilon_{n}}\right)^{-\frac{1}{s_{k}}}\right)$ as given by Lemma \ref{2.1HS}.
 Let $\vartheta_{k}$ denote a polar coordinate around $\bar{b}_{k}$, on
each $B_{\frac{T_{\varepsilon_{n}}}{10d_{k}}}\left(x_{j}\right)$, 
{according to notation \eqref{annulus}}, we define $U_{\varepsilon_{n}} (x) =U_{\varepsilon_{n}}^{j,k}(x)$ in $B_{\frac{T_{\varepsilon_{n}}}{10d_{k}}}\left(\bar{b}_{k}\right)$   where 
\begin{equation}\label{up1}
  U_{\varepsilon_{n}}^{j,k}(x)=
    \begin{cases}
    \frac{\left|x-x_{j}\right|}{\lambda\varepsilon}f_{0}(\lambda)e^{i\vartheta_{j}} & \text{on $B_{\lambda\varepsilon_{n}}\left(x_{j}\right)$} \\
    f_{0}\left(\frac{\left|x-x_{j}\right|}{\varepsilon_{n}}\right)e^{i\vartheta_{j}} &
                               \text {on $B_{\lambda\varepsilon_{n},\frac{T_{\varepsilon_{n}}}{20d_{k}}}\left(x_{j}\right)$}
  \\
  \left(f_{0}\left(\frac{T_{\varepsilon_{n}}}{20d_{k}\varepsilon_{n}}\right)+\left(\frac{|x-x_{j}|-\frac{T_{\varepsilon_{n}}}{20d_{k}}}{\frac{T_{\varepsilon_{n}}}{20d_{k}}}\right)\left(1-f_{0}\left(\frac{T_{\varepsilon_{n}}}{20d_{k}\varepsilon_{n}}\right)\right)\right)e^{i \vartheta_{j}}&\text{on $B_{\frac{T_{\varepsilon_{n}}}{20d_{k}},\frac{T_{\varepsilon_{n}}}{10d_{k}}}(\bar{b}_{k})$}.
    \end{cases}
\end{equation}
In this step we prove that
\begin{equation}\label{refmap1ter}
\begin{array}{ll}
  \displaystyle E_{\varepsilon_{n}}\left(U_{\varepsilon_{n}}^{j,k},B_{\frac{T_{\varepsilon_{n}}}{10d_{k}}}\left(x_{j}\right)\right)&    \displaystyle \leq -2\pi p_{0}\frac{1}{s_{k}}\log\log\frac{1}{\varepsilon_{n}}+2\pi p_{0}\log\frac{1}{\varepsilon_{n}}\\
  \\
&  \displaystyle -2\pi
    p_{0}I\left(\frac{1}{\varepsilon_{n}}\left(\log\frac{1}{\varepsilon_{n}}\right)^{-\frac{1}{s_{k}}}\right)+O(1).
\end{array}
\end{equation}
To this aim let us observe that of course we have
\begin{equation}\label{up2}
   E_{\varepsilon_{n}}\left(U_{\varepsilon_{n}}^{j,k},B_{\lambda\varepsilon_{n}}\left(x_{j}\right)\right)= O(1).
\end{equation}
By putting $U_{\varepsilon_{n}}^{j,k}(x)$ in the energy we obtain
\begin{equation}\label{a}
\begin{array}{ll}
    \displaystyle E_{\varepsilon_{n}}\left(U_{\varepsilon_{n}}^{j,k},B_{\lambda\varepsilon_{n},\frac{T_{\varepsilon_{n}}}{20d_{k}}}\left(x_{j}\right)\right)=&    \displaystyle 
2\pi\int_{\lambda\varepsilon_{n}}^{\frac{T_{\varepsilon_{n}}}{20d_{k}}}pf_{0}^{'2}rdr
    +\underbrace{2\pi\int_{\lambda\varepsilon_{n}}^{\frac{T_{\varepsilon_{n}}}{20d_{k}}}p\frac{f_{0}^{2}}{r}dr}_{(a)}+\\
    \\
  &  \displaystyle+\frac{2\pi}{\varepsilon^{2}}\int_{\lambda\varepsilon_{n}}^{\frac{T_{\varepsilon_{n}}}{20d_{k}}}J\left(1-f_{0}^{2}\right)rdr.    
\end{array}
\end{equation}
By Lemma \ref{2.5HS} and \eqref{HS3.2} we
deduce
\begin{equation}\label{up12}
   \int_{\lambda\varepsilon_{n}}^{\frac{T_{\varepsilon_{n}}}{20d_{k}}}pf_{0}^{'2}rdr\leq C
\end{equation}
and
\begin{equation}\label{up13}
    \frac{1}{\varepsilon^{2}}\int_{\lambda\varepsilon_{n}}^{\frac{T_{\varepsilon_{n}}}{20d_{k}}}J\left(1-f_{0}^{2}\right)rdr\leq C.
\end{equation}
Hence let us split term $(a)$ in \eqref{a} in the following way
\begin{equation}\label{up3}
   (a)=2\pi\int_{\lambda\varepsilon_{n}}^{\frac{T_{\varepsilon_{n}}}{20d_{k}}}p\frac{f_{0}^{2}}{r}dr=
   \underbrace{2\pi\int_{\lambda\varepsilon_{n}}^{\frac{T_{\varepsilon_{n}}}{20d_{k}}}(p-p_{0})\frac{f_{0}^{2}}{r}dr}_{(1)}+
   \underbrace{2\pi p_{0}\int_{\lambda\varepsilon_{n}}^{\frac{T_{\varepsilon_{n}}}{20d_{k}}}\frac{f_{0}^{2}}{r}dr}_{(2)}.
\end{equation}
Let us observe that
$$\left|x-\bar{b}_{k}\right|^{s_{k}}\leq 2^{s_{k}}\left(\left|x-x_{j}\right|^{s_{k}}+\left|x_{j}-\bar{b}_{k}\right|^{s_{k}}\right)\quad \forall j\in\{1,...,d_k\},$$
hence,  by $\left(\ref{p}\right)$  we have
\begin{equation*}
\begin{split}
    (1)\leq &\, \frac{2^{s_{k}+1}}{(20d_{k})^{s_{k}}}\pi\beta_{k}\left(\log\frac{1}{\varepsilon_{n}}\right)^{-1}\int_{\lambda\varepsilon}^{\frac{T_{\varepsilon_{n}}}{20d_{k}}}\frac{f_{0}^{2}}{r}dr+
 2\pi\beta_{k}\left(\log\frac{1}{\varepsilon_{n}}\right)^{-1}\int_{\lambda\varepsilon_{n}}^{\frac{T_{\varepsilon_{n}}}{20d_{k}}}\frac{f_{0}^{2}}{r}dr
\\
   = 
    &\,-2\pi\beta_{k}\left(\frac{1}{10^{s_{k}}d_{k}^{s_{k}}}+1\right)\left(\log\frac{1}{\varepsilon_{n}}\right)^{-1}\int_{\lambda\varepsilon_{n}}^{\frac{T_{\varepsilon_{n}}}{20d_{k}}}\frac{1-f_{0}^{2}}{r}dr+\\
    +&\,
 2\pi\beta_{k}\left(\frac{1}{10^{s_{k}}d_{k}^{s_{k}}}+1\right)\left(\log\frac{1}{\varepsilon_{n}}\right)^{-1}\int_{\lambda\varepsilon_{n}}^{\frac{T_{\varepsilon_{n}}}{20d_{k}}}\frac{dr}{r}.
  \end{split}
  \end{equation*}
By Lemma \ref{2.1HS} and Lemma \ref{2.7HS} 

  \begin{equation*}
  \begin{split}
(1)&\leq
    -2\pi\beta_{k}\left(\frac{1}{10^{s_{k}}d_{k}^{s_{k}}}+1\right)\left(\log\frac{1}{\varepsilon_{n}}\right)^{-1}I\left(\frac{1}{\varepsilon_{n}}\left(\log\frac{1}{\varepsilon_{n}}\right)^{-\frac{1}{s_{k}}}\right)
\\&
  \,\,\,\,\,\,+2\pi\beta_{k}\left(\frac{1}{10^{s_{k}}d_{k}^{s_{k}}}+1\right)\left(\log\frac{1}{\varepsilon_{n}}\right)^{-1}\left[-\frac{1}{s_{k}}\log\log\frac{1}{\varepsilon_{n}}+\log\frac{1}{\lambda\varepsilon_{n}}\right]+ O(1)\\
  &=-2\pi\beta_{k}\left(\frac{1}{10^{s_{k}}d_{k}^{s_{k}}}+1\right)\left(\log\frac{1}{\varepsilon_{n}}\right)^{-1}I\left(\frac{1}{\varepsilon_{n}}\left(\log\frac{1}{\varepsilon_{n}}\right)^{-\frac{1}{s_{k}}}\right)\\
  &\,\,\,-2\dfrac{\pi\beta_{k}}{s_{k}}\left(\frac{1}{10^{s_{k}}d_{k}^{s_{k}}}+1\right)\left(\log\frac{1}{\varepsilon_{n}}\right)^{-1}\log\log\frac{1}{\varepsilon_{n}} 
\\&
 \,\,\,\,\,\,+2\pi\beta_{k}\left(\frac{1}{10^{s_{k}}d_{k}^{s_{k}}}+1\right)\left(1-\log \lambda \left(\log\frac{1}{\varepsilon_{n}}\right)^{-1} \right)+O(1).
\end{split}
\end{equation*}
Let us observe that
\begin{equation*}
    \lim_{n\rightarrow
    +\infty}\left(\log\frac{1}{\varepsilon_{n}}\right)^{-1}\log\log\frac{1}{\varepsilon_{n}}=0
\end{equation*}
and again by $\left(\ref{HS2.12}\right)$ that \begin{equation*}
    \lim_{n\rightarrow
    +\infty}\left(\log\frac{1}{\varepsilon_{n}}\right)^{-1}I\left(\frac{1}{\varepsilon_{n}}\left(\log\frac{1}{\varepsilon_{n}}\right)^{-\frac{1}{s_{k}}}\right)=0.
\end{equation*}
Then we can conclude
\begin{equation*}
    (1)\leq O(1).
\end{equation*}
Now let us consider the second term in the right hand side of $\left(\ref{up3}\right)$
\begin{equation*}
\begin{split}
    (2)=&2\pi p_{0}\int_{\lambda\varepsilon_{n}}^{\frac{T_{\varepsilon_{n}}}{20d_{k}}}\frac{f_{0}^{2}}{r}dr=
    -2\pi p_{0}\int_{\lambda\varepsilon_{n}}^{\frac{T_{\varepsilon_{n}}}{20d_{k}}}\frac{1-f_{0}^{2}}{r}dr+
    2\pi p_{0}\int_{\lambda\varepsilon_{n}}^{\frac{T_{\varepsilon_{n}}}{20d_{k}}}\frac{dr}{r}\\=&
    -2\pi p_{0}I\left(\frac{1}{\varepsilon_{n}}\left(\log\frac{1}{\varepsilon_{n}}\right)^{-\frac{1}{s_{k}}}\right)+
     2\pi p_{0}\left(-\frac{1}{s_{k}}\log\log\frac{1}{\varepsilon_{n}}+\log\frac{1}{\lambda\varepsilon_{n}}\right)+
     O(1)\\=&
    -2\pi p_{0}I\left(\frac{1}{\varepsilon_{n}}\left(\log\frac{1}{\varepsilon_{n}}\right)^{-\frac{1}{s_{k}}}\right)
     -2\pi p_{0}\frac{1}{s_{k}}\log\log\frac{1}{\varepsilon_{n}}+2\pi p_{0}\log\frac{1}{\varepsilon_{n}}+
     O(1).
     \end{split}
\end{equation*}
By collecting together, we get
\begin{equation}\label{up14}
  (a)=(1)+(2)\leq
  -2\pi p_{0}I\left(\frac{1}{\varepsilon_{n}}\left(\log\frac{1}{\varepsilon_{n}}\right)^{-\frac{1}{s_{k}}}\right)-2\pi p_{0}\frac{1}{s_{k}}\log\log\frac{1}{\varepsilon_{n}}+2\pi p_{0}\log\frac{1}{\varepsilon_{n}}+
     O(1).
\end{equation}
Let us observe that $\left(\ref{refmap1ter}\right)$ will follows from $\left(\ref{up2}\right)$, $\left(\ref{up12}\right)$,
$\left(\ref{up13}\right)$ and
$\left(\ref{up14}\right)$ once we prove that
\begin{equation}\label{refmap1**}
E_{\varepsilon_{n}}\left(U_{\varepsilon_{n}}^{j,k},B_{\frac{T_{\varepsilon_{n}}}{20d_{k}},\frac{T_{\varepsilon_{n}}}{10d_{k}}}(x_{j})\right)\leq C.
\end{equation}
In order to verify $\left(\ref{refmap1**}\right)$ we write,
\begin{equation}
 U_{\varepsilon_{n}}^{j,k}(x_{j}+r e^{i\vartheta_{j}})=z(r)e^{i\vartheta_{j}}\,\,\,\, \text{on}\, B_{\frac{T_{\varepsilon_{n}}}{20d_{k}},\frac{T_{\varepsilon_{n}}}{10d_{k}}}(x_{j})
\end{equation}
where
$$
z(r)=f_{0}\left(\frac{T_{\varepsilon_{n}}}{20d_{k}\varepsilon_{n}}\right)+\left(\frac{r-\frac{T_{\varepsilon_{n}}}{20d_{k}}}{\frac{T_{\varepsilon_{n}}}{20d_{k}}}\right)\left(1-f_{0}\left(\frac{T_{\varepsilon_{n}}}{20d_{k}\varepsilon_{n}}\right)\right).
$$
Acting as in Proposition 3.1 in \cite{HS1}, by the properties of $f_{0}$ of Lemma \ref{2.1HS} and as $T_{\varepsilon_{n}}$ go to zero when $\varepsilon_{n}$ tends to zero, we compute
\begin{equation}\label{U1}
\begin{split}
\displaystyle \int_{B_{\frac{T_{\varepsilon_{n}}}{20d_{k}},\frac{T_{\varepsilon_{n}}}{10d_{k}}}(x_{j})}|\nabla U_{\varepsilon_{n}}^{j,k}|^{2}=&\displaystyle\int_{B_{\frac{T_{\varepsilon_{n}}}{20d_{k}},\frac{T_{\varepsilon_{n}}}{10d_{k}}}(x_{j})} z^{2}|\nabla \vartheta_{k}|^{2}+2\pi \displaystyle\int_{\frac{T_{\varepsilon_{n}}}{20d_{k}}}^{\frac{T_{\varepsilon_{n}}}{10d_{k}}}\left(z'\right)^{2}r dr\\
    =&
    O(1)+2\pi \left(\dfrac{1-f_{0}\left(\frac{T_{\varepsilon_{n}}}{20d_{k}\varepsilon_{n}}\right)}{\eta_{0}}\right)^{2}\displaystyle\int_{\frac{T_{\varepsilon_{n}}}{20d_{k}}}^{\frac{T_{\varepsilon_{n}}}{10d_{k}}} r dr\leq C.
\end{split}
\end{equation} 
About the second term of the energy, using the inequality $J(t)\leq tj(t)$, Lemma \ref{2.1HS} and Lemma \ref{2.2HS}, we obtain
\begin{equation}\label{U2}
\begin{split}
\dfrac{1}{\varepsilon_{n}^{2}}\displaystyle\int_{B_{\frac{T_{\varepsilon_{n}}}{20d_{k}},\frac{T_{\varepsilon_{n}}}{10d_{k}}}(x_{j})} J\left(1-|U_{\varepsilon_{n}}^{j,k}|^{2}\right)\leq & \dfrac{C}{\varepsilon_{n}^{2}}\displaystyle\int_{B_{\frac{T_{\varepsilon_{n}}}{20d_{k}},\frac{T_{\varepsilon_{n}}}{10d_{k}}}(x_{j})}j\left(1-|U_{\varepsilon_{n}}^{j,k}|^{2}\right)\\ \leq &
    \dfrac{C}{\varepsilon_{n}^{2}}j\left(1-f_{0}^{2}\left(\frac{T_{\varepsilon_{n}}}{20d_{k}\varepsilon_{n}}\right)\right)\left(\frac{T_{\varepsilon_{n}}^{2}}{100d_{k}^{2}}-\frac{T_{\varepsilon_{n}}^{2}}{400d_{k}^{2}}\right)
    \\ =&\dfrac{C}{\varepsilon_{n}^{2}}\dfrac{3}{\lambda\left(\frac{T_{\varepsilon_{n}}}{20d_{k}\varepsilon_{n}}\right)^{2}}\dfrac{T_{\varepsilon_{n}}^{2}}{400}=O(1).
\end{split}
\end{equation}
Hence by $(\ref{U1})$ and $(\ref{U2})$ we get $(\ref{refmap1**})$.\\
Finally, by $\left(\ref{up2}\right)$, $\left(\ref{up12}\right)$,
$\left(\ref{up13}\right)$,
$\left(\ref{up14}\right)$ and $(\ref{refmap1**})$ we can write 
\begin{equation}\label{up4}
  E_{\varepsilon_{n}}\left(U_{\varepsilon_{n}}^{j,k},B_{\frac{T_{\varepsilon_{n}}}{10d_{k}}}\left(x_{j}\right)\right)\leq
      -2\pi p_{0}\frac{1}{s_{k}}\log\log\frac{1}{\varepsilon_{n}}+2\pi
  p_{0}\log\frac{1}{\varepsilon_{n}}-2\pi p_{0}I\left(\frac{1}{\varepsilon_{n}}\left(\log\frac{1}{\varepsilon_{n}}\right)^{-\frac{1}{s_{k}}}\right)+O(1).
\end{equation}
\textbf{Step 3.} We define the function
$U_{\varepsilon_{n}}$  in $\bigcup_{j=1}^{d_{k}}B_{T_{\varepsilon_{n}}}\left(x_{j}\right)$ such that
\begin{equation*}
    U_{\varepsilon_{n}}^{k}(x)=U_{\varepsilon_{n}}^{j,k}(x) \quad \hbox{if} \quad x\in B_{T_{\varepsilon_{n}}}\left(x_{j}\right).
\end{equation*}
As the discs centered in $x_{j}$ are disjoint and as they are exactly $d_{k}$ discs  we get
\begin{equation}\label{up15}
\begin{array}{ll}
    \displaystyle E\left(U_{\varepsilon_{n}}^{k}(x),\bigcup_{j=1
}^{d_{k}}B_{T_{\varepsilon_{n}}}\left(x_{j}\right)\right)\leq &\displaystyle
-2\pi
p_{0}d_{k}I\left(\frac{1}{\varepsilon_{n}}\left(\log\frac{1}{\varepsilon_{n}}\right)^{-\frac{1}{s_{k}}}\right)
  -2\pi p_{0}\frac{d_{k}}{s_{k}}\log\log\frac{1}{\varepsilon_{n}}\\
  &\displaystyle+2\pi p_{0}d_{k}\log\frac{1}{\varepsilon_{n}}+O(1).
\end{array}
\end{equation}
By $\left(\ref{refmap1}\right)$, $\left(\ref{refmap2}\right)$ and
$\left(\ref{up15}\right)$ we have
\begin{equation}\label{up10}
\begin{split}
  E_{\varepsilon_{n}}\left(U_{\varepsilon_{n}}^{k},B_{\eta_{0}}(\bar{b}_{k})\right)\leq & 2\pi
  p_{0}\dfrac{d_{k}^{2}}{s_{k}}\log\log\dfrac{1}{\varepsilon_{n}}-2\pi p_{0}d_{k}I\left(\frac{1}{\varepsilon_{n}}\left(\log\frac{1}{\varepsilon_{n}}\right)^{-\frac{1}{s_{k}}}\right)\\
&
  -2\pi p_{0}\dfrac{d_{k}}{s_{k}}\log\log\dfrac{1}{\varepsilon_{n}}+2\pi
  p_{0}d_{k}\log\dfrac{1}{\varepsilon_{n}}+O(1).
\end{split}
\end{equation}
Finally, we pose $U_{\varepsilon_{n}}(x)=w$ on $G\setminus
\bigcup_{k=1}^{N}\overline{B_{\eta_{0}}(\bar{b}_{k})}$ where $w$ is any
$S^{1}$-valued map of class $C^{1}$ on this domain which equals $g$
on $\partial G$ and
$\left(\frac{x-\bar{b}_{k}}{\left|x-\bar{b}_{k}\right|}\right)^{d_{k}}$
on $\partial B_{\eta_{0}}(\bar{b}_{k})$ for $k=1,.., N$. Then $U_{\varepsilon_{n}}\in H^{1}_{g}\left(G,\textbf{C}\right)$ and we get
\begin{equation}\label{up11}
\begin{split}
  E_{\varepsilon_{n}}\left(u_{\varepsilon_{n}}\right)\leq E_{\varepsilon_{n}}\left(U_{\varepsilon_{n}}\right)\leq & 2\pi p_{0}d\log\frac{1}{\varepsilon_{n}}+2\pi
    p_{0}\Sigma_{k=1}^{N}\frac{d_{k}^{2}-d_{k}}{s_{k}}\log\log\frac{1}{\varepsilon_{n}}\\
 &
-2\pi
p_{0}dI\left(\frac{1}{\varepsilon_{n}}\left(\log\frac{1}{\varepsilon_{n}}\right)^{-\frac{1}{s_{k}}}\right)+O(1)
\end{split}
\end{equation}
which is $\left(\ref{up11***}\right)$.
\subsection{A lower bound for the energy}\label{sub3.2}
When $G$ is starshaped, using a Pohozaev identity, we obtain
\begin{equation}\label{HS3.2}
\dfrac{1}{\varepsilon^{2}}\int_{G}J\left(1-\left| u_{\varepsilon}\right|^{2}\right)\leq C_{0},\,\, \forall \varepsilon >0.
\end{equation}

By following the same arguments of Lemma 3.1 and Lemma 3.2 in \cite{HS1} we get
\begin{equation}\label{grad}
\left\|u_{\varepsilon}\right\|_{L^{\infty}(G)}\leq 1 \,\text{and}\,\,\left\|\nabla u_{\varepsilon}\right\|_{L^{\infty}(G)}\leq \dfrac{C}{\varepsilon}.
\end{equation}  
Using the construction in \cite{BBH2} we know that there
exist $\lambda
>0$ and a collection of balls $\left\lbrace B_{\lambda
\varepsilon}\left( y_{j}^{\varepsilon }\right)\right\rbrace_{j\in J}$ such that
\begin{equation}\label{bd}
\left\{ x\in \overline{G}:\left\vert u_{\varepsilon }\left( x\right)
\right\vert \leq \frac{3}{4}\right\} \subset
\bigcup_{j\in J}B_{\lambda
\varepsilon}\left( y_{j}^{\varepsilon }\right),
\end{equation}
\begin{equation*}
\left|y_{i}^{\varepsilon}-y_{j}^{\varepsilon}\right|\geq 8\lambda \varepsilon \,\,\forall i,j\in J, i\neq j
\end{equation*}
and
$$\text{card } J\leq N_{b}.$$
By construction, the degrees 
$$\nu_{j}=\deg \left( u_{\varepsilon
},\partial B_{\lambda \varepsilon}\left( y_{j}^{\varepsilon} \right)\right), j\in J$$ 
are well defined.
Given any subsequence $\varepsilon_{n}\rightarrow 0$ we may extract a subsequence (still denoted by $\varepsilon_{n}$) such that
$$\text{card } J_{\varepsilon_{n}}=cost=N_{1}$$
and
\begin{equation}\label{defyi} y_{j}=y_{j}^{\varepsilon _{n}}\rightarrow l_{j}\in \overline{G},\,j=1,..,N_{1}.
\end{equation}  
Let $\underline{b}_{1}, \underline{b}_{2},..,\underline{b}_{N_{2}}$ be the distinct points among the $\left\lbrace l_{j}\right\rbrace_{j=1}^{N_{1}}$ and set $$I_{k}=\left\lbrace j\in \left\lbrace 1,..,N_{1}\right\rbrace; y_{j}^{\varepsilon_{n}}\rightarrow \underline{b}_{k}\right\rbrace, k=1,.., N_{2}.$$
Denoting by $d_{k}=\sum_{j\in I_{k}}\nu_{j}$ for every $k=1,..,N_{2}$, we clearly have and $\sum_{k=1}^{N_{2}}d_{k}=d$. By following the same arguments as in \cite{AS}, thanks to the previous upper bound, applied to $\bar{b}_k=\underline{b}_k$, and Proposition  \ref{3.3HS}, we get 
\begin{equation}\label{dk}
d_{k}>0\,\, \text{for every}\, k=1,..,N_{2}
\end{equation}
and
\begin{equation}
\underline{b}_{k}\in \Lambda=p^{-1} (p_0)\,\, \text{for every}\,  k=1,..,N_{2}.
\end{equation} 
Hence,  having in mind \eqref{P*}, in the following we can set $N_{2}=N$ and $\underline{b}_k=b_k$. Moreover acting as in \cite{AS}, Lemma 2.1 by Proposition \ref{3.2HS} and Proposition \ref{3.3HS}, we get $\nu_{j}=+1$ for every $j\in I_{k}$. \\
Let $\eta$ satisfy
\begin{equation}\label{eta}
    0<\eta<\frac{1}{2}{\min }\left(\underset{i\neq j}{\min }\left|b_{i}-b_{j}\right|,\underset{i=1,..,N_{2}}{\min }dist\left(b_{i},\partial
    G\right)\right).
\end{equation}
and take 
$T_{\varepsilon_{n}}$ as in $\eqref{definitionTepsilon}$. 
We now are able to prove the following lower bound : 
\begin{proposition}\label{proplb} Assume $G$ is starshaped and $\eqref{P*}$ and $\eqref{p}$ hold true.  
Then we have, for a subsequence $\varepsilon_{n}\rightarrow 0$ 
\begin{equation}\label{lb}
    \begin{split}
   E_{\varepsilon_{n}}\left(u_{\varepsilon_{n}}\right)\geq &\, 2\pi
p_{0}d\log\frac{1}{\varepsilon_{n}}+2\pi
p_{0}\sum_{k=1}^{N}\frac{d_{k}^{2}-d_{k}}{s_{k}}\log\log\frac{1}{\varepsilon_{n}} 
-2\pi
p_{0}dI\left(\frac{1}{\varepsilon_{n}}\left(\log\frac{1}{\varepsilon_{n}}\right)^{-\frac{1}{s_{k}}}\right) \\
&
-2\pi
    p_{0}\sum_{k=1}^{N}d_{k}^{2}I\left(\left(\log\frac{1}{\varepsilon_{n}}\right)^{\frac{1}{s_{k}}}\right)
    + {9\over 8}\pi p_{0} \sum_{k=1}^{N}\sum_{i\neq j}\log \dfrac{ \left(\log\frac{1}{\varepsilon_{n}}\right)^{-\frac{1}{s_{k}}} }{\mid y_{i}-y_{j}\mid}+O(1),
\end{split}
\end{equation}
where the points $y_i$ and $y_j$, $i,\,j \in I_k$,  are as in \eqref{defyi}.
\end{proposition}

\textbf{Proof.} The proof developes into two steps.\\
\textbf{Step1.} {By following a similar argument as in \cite{AS},} at first we prove 
\begin{equation}\label{AS3.13}
    \max_{i\in I_{k}}\left|b_{k}-y_{i}\right|=R_n\sim |\log\varepsilon_{n}|^{-\frac{1}{s_{k}}}
\end{equation}
for every 
$k=1,..,N_{2}$ with $\left|I_{k}\right| = d_{k}>1$.\\
We know that $B_{\eta}\left(b_{k}\right)$ contains
exactly $d_{k}$ bad discs $B_{\lambda\varepsilon_{n}}\left(y_{i}\right)$,
such that for every $\alpha\in\left(0,1\right)$
\begin{equation}\label{epsalfa}
\left|y_{i}-y_{j}\right|> \varepsilon_{n}^{\alpha}\quad \forall i\neq j.
\end{equation}
For any fixed $\alpha\in\left(0,1\right)$, we have
\begin{equation}\label{3.1}
\begin{split}
E\left(u_{\varepsilon_{n}},B_{\eta}\left(b_{k}\right)\right)
  \geq & E\left(u_{\varepsilon_{n}}, B_{2R_{n},\eta}\left(b_{k}\right)\right)+ E\left(u_{\varepsilon},B_{2R_{n}}\left(b_{k}\right)\setminus\bigcup_{i\in I_{k}} B_{\varepsilon_{n}^{\alpha}}\left(y_{i}\right)\right)
  \\
   &+E\left(u_{\varepsilon_{n}},\bigcup_{i\in I_{k}} B_{\lambda\varepsilon_{n},\varepsilon_{n}^{\alpha}} \left( y_{i}\right)\right)
    = (a) + (b)+ (c).
\end{split}
\end{equation}
Taking into account \eqref{HS3.2}, by Proposition \ref{3.2HS}, there exist two constants $C_{1}$ and $C_{3}$ depending only on $C_{0}$ and a constant $C_{2}$ depending on $C_{0}$ and $d_{k}$, such that
\begin{equation}\label{3.2}
   (a)\geq 2\pi
    d_{k}^{2}p_{0}\left[\log\frac{\eta}{2R_{n}}-I\left(\frac{\eta}{2R_{n}}\right)\right]-d_{k}^{2}C_{1},
\end{equation}
\begin{equation}\label{3.3}
    (b)\geq 2\pi
    d_{k}p_{0}\left[\log\frac{2R_{n}}{\varepsilon_n^{\alpha}}
    -I\left(\frac{2R_{n}}{\varepsilon_{n}^{\alpha}}\right)\right]-C_{2}
\end{equation}
and
\begin{equation}\label{3.4}
  (c)\geq 2\pi
    \left(d_{k}-1\right)p_{0}\left[\log\frac{\varepsilon_n^{\alpha}}{\lambda\varepsilon_{n}}
    -I\left(\frac{\varepsilon_{n}^{\alpha}}{\lambda\varepsilon_{n}}\right)\right]+ 
  2\pi
    \left(p_{0}+\alpha_{k}\frac{R_{n}^{s_{k}}}{4}\right)\left[\log\frac{\varepsilon_{n}^{\alpha}}{\lambda\varepsilon_{n}}
    -I\left(\frac{\varepsilon_{n}^{\alpha}}{\lambda\varepsilon_{n}}\right)\right]-C_{3}.
    \end{equation}
Let us denote 
\begin{equation}\label{deff}
    f\left(R_{n}\right)=2\pi p_{0}d_{k}\log\frac{1}{\varepsilon_{n}}+2\pi
    p_{0}\left(d_{k}^{2}-d_{k}\right)\log\frac{1}{R_{n}}+\frac{\pi}{2}\alpha_{k}\left(1-\alpha\right)R_{n}^{s_{k}}\log\frac{1}{\varepsilon_{n}}
\end{equation}
and
\begin{equation}\label{defg}
    g\left(R_{n}\right)=2\pi
    d_{k}^{2}p_{0}I\left(\frac{1}{R_{n}}\right)+2\pi
    d_{k}p_{0}I\left(\frac{R_{n}}{\varepsilon_{n}^{\alpha}}\right)+2\pi\left(p_{0}+\alpha_{k}\frac{R_{n}^{s_{k}}}{4}\right)I\left(\frac{1}{\varepsilon_{n}^{1-\alpha}}\right)+C_{4}.
\end{equation}
where $C_{4}$ is a constant depending only on $C_{0}$ and $d_{k}$.
Then
\begin{equation}\label{3.5}
    E\left(u_{\varepsilon_{n}},B_{\eta}\left(b_{k}\right)\right)\geq
    f\left(R_{n}\right)-g\left(R_{n}\right)-C_{4}.
\end{equation}
Now let us observe that, for for $n$ large enough, we get 
$$
  \dfrac{\eta}{2R_{n}}\geq 1,
$$
since $R_n$ tends to $0$. Moreover by \eqref{epsalfa} it holds
$$\varepsilon_{n}^{\alpha} < \left|y_{i}-y_{j}\right| \leq \left|y_i-b_{k}\right| + \left|y_j-b_{k}\right| \leq 2R_n\quad \forall i\not= j.$$ 
Hence we get
\begin{equation}\label{3.6}
\dfrac{\varepsilon_{n}^{\alpha}}{2}\leq R_{n}\leq \frac{\eta}{2} 
\end{equation}


Let us pose $R_{n}=c_{n}\left(\log\dfrac{1}{\varepsilon_{n}}\right)^{-\frac{1}{s_{k}}}$ and consider the following difference 
\begin{equation*}
\begin{array}{c}
  \left[f\left(R_{n}\right)-g\left(R_{n}\right)\right]-
    \left[f\left(\left(\log\dfrac{1}{\varepsilon_{n}}\right)^{-\frac{1}{s_{k}}}\right)
    -g\left(\left(\log\dfrac{1}{\varepsilon_{n}}\right)^{-\frac{1}{s_{k}}}\right)\right]= \\
    \\
  \underbrace{\left[f\left(R_{n}\right)-f\left(\left(\log\frac{1}{\varepsilon_{n}}\right)^{-\frac{1}{s_{k}}}\right)\right]}_{(1)}+
 \underbrace{\left[g\left(\left(\log\frac{1}{\varepsilon_{n}}\right)^{-\frac{1}{s_{k}}}\right)-g\left(R_{n}\right)\right]}_{(2)}.
\end{array}
\end{equation*}
By \eqref{deff} and \eqref{defg} we get
\begin{equation}\label{3.6*}
    (1)=2\pi p_{0}\left(d_{k}^{2}-d_{k}\right)\log\frac{1}{c_{n}}+\frac{\pi}{2}\alpha_{k}\left(1-\alpha\right)\left(c_{n}^{s_{k}}-1\right)
\end{equation}
and
\begin{equation}\label{3.7}
\begin{array}{ll}
  (2)=& 
  \displaystyle2\pi d_{k}^{2}p_{0}\left(I\left(\left(\log\dfrac{1}{\varepsilon_{n}}\right)^{\frac{1}{s_{k}}}\right)
  -I\left(\dfrac{1}{c_{n}}\left(\log\dfrac{1}{\varepsilon_{n}}\right)^{\frac{1}{s_{k}}}\right)\right)+\\
  \\
 & \displaystyle+2\pi p_{0}d_{k}\left(I\left(\dfrac{1}{\varepsilon_{n}^{\alpha}}\left(\log\dfrac{1}{\varepsilon_{n}}\right)^{-\frac{1}{s_{k}}}\right)
  -I\left(\dfrac{c_{n}}{\varepsilon_{n}^{\alpha}}\left(\log\dfrac{1}{\varepsilon_{n}}\right)^{-\frac{1}{s_{k}}}\right)\right)+ \\
  \\
  &\displaystyle
  +\dfrac{\pi\alpha_{k}\left(1-c_{n}^{s_{k}}\right)}{2}\left(\log\dfrac{1}{\varepsilon_{n}}\right)^{-1}
  I\left(\dfrac{1}{\varepsilon_{n}^{1-\alpha}}\right).\\
  \end{array}
    \end{equation}
Let us consider the case $c_{n}>1$. Therefore we have
\begin{equation}\label{R1}
R_{n}>\left(\log\frac{1}{\varepsilon_{n}}\right)^{-\frac{1}{_{s_{k}}}}.
\end{equation}

By $(\ref{HS2.11})$, $(\ref{R1})$ and as the functions $j^{-1}$ and $I$ are increasing, we get 
\begin{equation*}
    \begin{split}
      (2)\geq &
       -2\pi p_{0}d_{k}\displaystyle\int_{\frac{\varepsilon_{n}^{2\alpha}}{R_{n}^{2}}}^{\left(\log\frac{1}{\varepsilon_{n}}\right)^{\frac{2}{s_{k}}}\varepsilon_{n}^{2\alpha}}\frac{j^{-1}(t)}{t}dt+\dfrac{\pi\alpha_{k}\left(1-c_{n}^{s_{k}}\right)}{2}\left(\log\dfrac{1}{\varepsilon_{n}}\right)^{-1}
  I\left(\dfrac{1}{\varepsilon_{n}^{1-\alpha}}\right)\\ \geq 
      &-2\pi
p_{0}d_{k} j^{-1}\left(\left(\log\frac{1}{\varepsilon_{n}}\right)^{\frac{2}{s_{k}}}\varepsilon_{n}^{2\alpha}\right)\log c_{n}^{2}+\dfrac{\pi\alpha_{k}\left(1-c_{n}^{s_{k}}\right)}{2}\left(\log\dfrac{1}{\varepsilon_{n}}\right)^{-1}
  I\left(\dfrac{1}{\varepsilon_{n}^{1-\alpha}}\right).
\end{split}
\end{equation*}
Since
$$\lim_{n} \left(\log\frac{1}{\varepsilon_{n}}\right)^{\frac{2}{s_{k}}}\varepsilon_{n}^{2\alpha}=0
$$
and by $(\ref{HS2.12})$
\begin{equation}\label{limit0}
\lim_{n} \left(\log\dfrac{1}{\varepsilon_{n}}\right)^{-1}
  I\left(\dfrac{1}{\varepsilon_{n}^{1-\alpha}}\right)=0,
\end{equation}
by regularity of function $j^{-1}$ and as  $j^{-1}(0)=0$, there exists $n_{0} $  such that for $n\geq n_{0}$ we have
  \begin{equation}\label{3.7*}
   (2)\geq  2\delta\pi
    p_{0}d_{k}\log\frac{1}{c_{n}}+\frac{\pi}{2}\alpha_{k}\left(1-c_{n}^{s_{k}}\right)\gamma.
\end{equation}
Then, by denoting
\begin{equation}\label{defh}
    h\left(R_{n}\right)= f\left(R_{n}\right)- g\left(R_{n}\right),
    \end{equation}
by (\ref{3.6*}) and (\ref{3.7*}) and choosing $\delta=\dfrac{1}{2}$ and $\gamma= \dfrac{1-\alpha}{2}$, we get
\begin{equation*}
 h\left(R_{n}\right)-h\left(\left(\log\frac{1}{\varepsilon_{n}}\right)^{-\frac{1}{s_{k}}}\right)\geq  2\pi p_{0}\left(d_{k}^{2}-\dfrac{d_{k}}{2}\right)\log\frac{1}{c_{n}}+ \frac{\pi}{8}\alpha_{k}\left(1-\alpha\right)\left(c_{n}^{s_{k}}-1\right).
\end{equation*}
Hence we get
\begin{equation}\label{3.7b}
h\left(R_{n}\right)-h\left(\left(\log\frac{1}{\varepsilon_{n}}\right)^{-\frac{1}{s_{k}}}\right)\rightarrow
    +\infty \,\text{as}\,\, c_{n}\rightarrow
    +\infty.
\end{equation}
Now let us suppose there exists a subsequence $(c_{n_k})_k$, still denoted by $(c_{n})$, such that    $c_{n}<1$. Up to a subsequence we have
\begin{equation}\label{R1*}
R_{n}<\left(\log\frac{1}{\varepsilon_{n}}\right)^{-\frac{1}{_{s_{k}}}}.
\end{equation}

By $(\ref{HS2.11})$, $(\ref{R1*})$ and as the functions $j^{-1}$ and $I$ are increasing, we get
\begin{equation*}
    \begin{split}
      (2)\geq &  -2\pi
p_{0}d_{k}^{2}\displaystyle\int_{R_{n}^{2}}^{\left(\log\frac{1}{\varepsilon_{n}}\right)^{-\frac{2}{s_{k}}}}\frac{j^{-1}(t)}{t}dt+\dfrac{\pi}{2}\alpha_{k}\left(1-c_{n}^{s_{k}}\right)\left(\log\dfrac{1}{\varepsilon_{n}}\right)^{-1}
  I\left(\dfrac{1}{\varepsilon_{n}^{1-\alpha}}\right)\\ \geq & 
      -2\pi
p_{0}d_{k}^{2} j^{-1}\left(\left(\log\frac{1}{\varepsilon_{n}}\right)^{-\frac{2}{s_{k}}}\right)\log \frac{1}{c_{n}^{2}}+\dfrac{\pi}{2}\alpha_{k}\left(1-c_{n}^{s_{k}}\right)\left(\log\dfrac{1}{\varepsilon_{n}}\right)^{-1}
  I\left(\dfrac{1}{\varepsilon_{n}^{1-\alpha}}\right).\end{split}
\end{equation*}
Since
$$\lim_{n} \left(\log\frac{1}{\varepsilon_{n}}\right)^{-\frac{2}{s_{k}}}=0
$$
and by $(\ref{HS2.12})$
$$
\lim_{n} \left(\log\dfrac{1}{\varepsilon_{n}}\right)^{-1}
  I\left(\dfrac{1}{\varepsilon_{n}^{1-\alpha}}\right)=0,
$$
similarly to the previous case, by regularity of function $j^{-1}$ and as  $j^{-1}(0)=0$ there exists $n_{0} $  such that for $n\geq n_{0}$ we have 
  \begin{equation}\label{3.7**}
   (2)\geq - 2\delta\pi
    p_{0}d_{k}^{2}\log\frac{1}{c_{n}}+\frac{\pi}{2}\alpha_{k}\left(1-c_{n}^{s_{k}}\right)\gamma \geq - 2\delta\pi
    p_{0}d_{k}^{2}\log\frac{1}{c_{n}}.
\end{equation}
Then, by denoting
\begin{equation*}
    h\left(R_{n}\right)= f\left(R_{n}\right)- g\left(R_{n}\right),
    \end{equation*}
by $(\ref{3.6*})$ and $(\ref{3.7**})$  we get
\begin{equation*}
 h\left(R_{n}\right)-h\left(\left(\log\frac{1}{\varepsilon_{n}}\right)^{-\frac{1}{s_{k}}}\right)\geq 2\pi p_{0}\left(d_{k}^{2}
 -d_{k}-\delta d_{k}^{2}\right)\log\frac{1}{c_{n}}+ \frac{\pi}{2}\alpha_{k}\left(1-\alpha\right)\left(c_{n}^{s_{k}}-1\right). 
\end{equation*}
Let us choose $\delta > 0$ such that $d_{k}^{2}
 -d_{k}-\delta d_{k}^{2}>1$ or equivalently $\delta< 1-\dfrac{1+d_{k}}{d_{k}^{2}}$. This is possible as $d_{k}>1$ and then $ 1-\dfrac{1+d_{k}}{d_{k}^{2}} > 0$. For this choice it holds

\begin{equation}\label{3.7b*}
h\left(R_{n}\right)-h\left(\left(\log\frac{1}{\varepsilon_{n}}\right)^{-\frac{1}{s_{k}}}\right)\rightarrow
    +\infty \,\text{as}\, \,\frac{1}{c_{n}}\rightarrow
    +\infty.
\end{equation}
By \eqref{3.7b} and \eqref{3.7b*}, in both cases we can conclude as in $\cite{AS}$
\begin{equation}\label{3.7bb}
h\left(R_{n}\right)-h\left(\left(\log\frac{1}{\varepsilon_{n}}\right)^{-\frac{1}{s_{k}}}\right)\rightarrow
    +\infty \,\text{as}\, \max\left(c_{n},\frac{1}{c_{n}}\right)\rightarrow
    +\infty.
\end{equation}
By $(\ref{3.5})$ we get
\begin{equation*}
h\left(R_{n}\right)-h\left(\left(\log\frac{1}{\varepsilon_{n}}\right)^{-\frac{1}{s_{k}}}\right)\leq E\left(u_{\varepsilon_{n}},B_{\eta}\left(b_{k}\right)\right)+C_4-h\left(\left(\log\frac{1}{\varepsilon_{n}}\right)^{-\frac{1}{s_{k}}}\right).
\end{equation*}

 We know that $\overline{b}_{k}=b_{j}$ for some $j\in \{ 1,..,N \} $. Hence by using the upper bound \eqref{up11***} of Proposition \ref{upb}, taking into account \eqref{deff}, \eqref{defg} and \eqref{defh}, since $\alpha<1$, we obtain
\begin{equation}\label{H}
\begin{array}{c}
 h\left(R_{n}\right)-h\left(\left(\log\dfrac{1}{\varepsilon_{n}}\right)^{-\frac{1}{s_{k}}}\right)\leq -
       2\pi p_{0}d_{k}I\left(\dfrac{1}{\varepsilon_{n}}\left(\log\dfrac{1}{\varepsilon_{n}}\right)^{-\frac{1}{s_{k}}}\right)+2\pi p_{0} d_{k}^{2}I\left(\left(\log\dfrac{1}{\varepsilon_{n}}\right)^{\frac{1}{s_{k}}}\right)+\\
       2\pi p_{0}d_{k}I\left(\dfrac{1}{\varepsilon_{n}^{\alpha}}\left(\log\dfrac{1}{\varepsilon_{n}}\right)^{-\frac{1}{s_{k}}}\right)+2\pi\left(p_{0}d_{k}+\dfrac{\alpha_{k}}{4}\left(\log\dfrac{1}{\varepsilon_{n}}\right)^{-1}\right)I\left(\dfrac{1}{\varepsilon_{n}^{1-\alpha}}\right)+O(1).
\end{array}
\end{equation}

\noindent By assumption $(H2)$ and \eqref{HS2.11} in Lemma \ref{2.4HS}, we deduce that the functional $I$ is increasing,  
thus for $n$ large enough, we get
$$
I\left(\frac{1}{\varepsilon_{n}}\left(\log\frac{1}{\varepsilon_{n}}\right)^{-\frac{1}{s_{k}}}\right)\geq I\left(\frac{1}{\varepsilon_{n}^{\alpha}}\left(\log\frac{1}{\varepsilon_{n}}\right)^{-\frac{1}{s_{k}}}\right),
$$
$$
I\left(\frac{1}{\varepsilon_{n}}\left(\log\frac{1}{\varepsilon_{n}}\right)^{-\frac{1}{s_{k}}}\right)\geq I\left(\left(\log\frac{1}{\varepsilon_{n}}\right)^{\frac{1}{s_{k}}}\right)
$$
and
$$
I\left(\frac{1}{\varepsilon_{n}}\left(\log\frac{1}{\varepsilon_{n}}\right)^{-\frac{1}{s_{k}}}\right)\geq I\left(\frac{1}{\varepsilon_{n}^{1-\alpha}}\right).
$$
Hence, by \eqref{limit0},  the leading term of the second member in (\ref{H}) is the negative one and we can conclude that
\begin{equation}\label{estimate}
    h\left(R_{n}\right)-h\left(\left(\log\frac{1}{\varepsilon_{n}}\right)^{-\frac{1}{s_{k}}}\right)\rightarrow
    -\infty \,\text{as}\, n\rightarrow
    +\infty.
\end{equation}
This is a contradiction with \eqref{3.7bb} 
and arguing as in $\cite{AS}$, $(\ref{3.7bb})$ directly implies (\ref{AS3.13}).

\textbf{Step 2.} Let $\eta$ as in \eqref{eta} and $T_{\varepsilon_{n}}$ as in \eqref{definitionTepsilon}. We know
that $B_{\eta}\left(b_{k}\right)$ contains exactly $d_{k}$ bad discs
$B_{\lambda\varepsilon}\left(y_{j}\right)$, $j\in I_{k}$ satisfying
$\left(\ref{AS3.13}\right)$. 

We have
\begin{equation}\label{up30}
\begin{split}
    E_{\varepsilon_{n}}\left(u_{\varepsilon_{n}},B_{\eta}\left(b_{k}\right)\right)\geq &
    E_{\varepsilon_{n}}\left(u_{\varepsilon_{n}},B_{\eta}\left(b_{k}\right)\setminus
    B_{T_{\varepsilon_{n}}}\left(b_{k}\right)\right)+\sum_{j\in
I_{k}}E_{\varepsilon_{n}}\left(u_{\varepsilon_{n}},B_{T_{\varepsilon_{n}}}\left(b_{k}\right)\setminus
    B_{\lambda\varepsilon_{n}}\left(y_{j}\right)\right)
    \\
=& E_{1}+E_{2}.
    \end{split}
\end{equation}
By Proposition \ref{3.2HS}, we have
\begin{equation*}
    E_{1}\geq 2\pi p_{0}d_{k}^{2}\log\frac{\eta}{T_{\varepsilon_{n}}}-2\pi
    p_{0}d_{k}^{2}I\left(\frac{\eta}{T_{\varepsilon_{n}}}\right)-d_{k}C_{6}.
\end{equation*}
where $C_{6}$ is a constant depending only on $C_{0}$.\\
Then 
\begin{equation}\label{lb1}
    E_{1}\geq 2\pi p_{0}\frac{d_{k}^{2}}{s_{k}}\log\log\frac{1}{\varepsilon_{n}}-2\pi
    p_{0}d_{k}^{2}I\left(\left(\log\frac{1}{\varepsilon_{n}}\right)^{\frac{1}{s_{k}}}\right)+O(1).
\end{equation}
By \eqref{BMR5*positif} in Remark \ref{positif} applied to $y_{1},...,y_{d_{k}}$, as  $\nu_{j}=\deg \left( u_{\varepsilon
},\partial B\left( y_{j},\lambda \varepsilon) \right)\right)=+1$ for every $j=1,..,d_{k}$ and by \eqref{bd}, we have
\begin{equation*}
     E_{2}\geq 2\pi p_{0} d_{k}\left( \log \dfrac{T_{\varepsilon_{n}}}{\lambda\varepsilon_{n}} -I\left( \dfrac{T_{\varepsilon_{n}}}{\lambda\varepsilon_{n}}\right)\right)+{9\over 8}\pi p_{0}\sum_{i\neq j}\log \dfrac{T_{\varepsilon_{n}}}{\mid y_{i}-y_{j}\mid}-C_{7}
\end{equation*}
where $C_{7}$ is a constant depending only on $d_{k}$, $C_{0}$, and $p_{0}$ where $C_{0}$ is introduce in (\ref{HS3.2}). 
Then
\begin{equation}\label{lb3}
\begin{split}
E_{2}\geq & -2\pi
p_{0}\dfrac{d_{k}}{s_{k}}\log\log\dfrac{1}{\varepsilon_{n}}+ 2\pi
p_{0}d_{k}\log\dfrac{1}{\varepsilon_{n}}
-2\pi
p_{0}d_{k}I\left(\dfrac{1}{\varepsilon_{n}}\left(\log\dfrac{1}{\varepsilon_{n}}\right)^{-\frac{1}{s_{k}}}\right)\\
&+{9\over 8}\pi p_{0}\sum_{i\neq j}\log \dfrac{T_{\varepsilon_{n}}}{\mid y_{i}-y_{j}\mid}+O(1).
 \end{split}
\end{equation}
By collecting together $\left(\ref{lb1}\right)$ and $\left(\ref{lb3}\right)$ we obtain
\begin{equation}\label{lbk} 
\begin{split}
    E_{\varepsilon_{n}}\left(u_{\varepsilon_{n}},B_{\eta}\left(b_{k}\right)\right)\geq & 2\pi p_{0}\dfrac{d_{k}^{2}-d_{k}}{s_{k}}\log\log\dfrac{1}{\varepsilon_{n}}-2\pi
    p_{0}d_{k}^{2}I\left(\left(\log\dfrac{1}{\varepsilon_{n}}\right)^{\frac{1}{s_{k}}}\right)+2\pi p_{0}d_{k}\log\dfrac{1}{\varepsilon_{n}}
    \\
     &-2\pi
p_{0}d_{k}I\left(\dfrac{1}{\varepsilon_{n}}\left(\log\dfrac{1}{\varepsilon_{n}}\right)^{-\frac{1}{s_{k}}}\right)+{9\over 8}\pi p_{0}\sum_{i\neq j}\log \dfrac{T_{\varepsilon_{n}}}{\mid y_{i}-y_{j}\mid}+O(1).
\end{split}
\end{equation}
Summing over $k$ we have
\begin{equation}\label{lb4}
\begin{split}
   E_{\varepsilon_{n}}\left(u_{\varepsilon_{n}}\right)\geq & E_{\varepsilon_{n}}\left(u_{\varepsilon_{n}},\bigcup_{k=1}^{N}B_{\eta}\left(b_{k}\right)\right)\geq 2\pi
p_{0}d\log\dfrac{1}{\varepsilon_{n}}
+2\pi
p_{0}\sum_{k=1}^{N}\dfrac{d_{k}^{2}-d_{k}}{s_{k}}\log\log\dfrac{1}{\varepsilon_{n}}\\
&-2\pi
    p_{0}\sum_{k=1}^{N}d_{k}^{2}I\left(\left(\log\dfrac{1}{\varepsilon_{n}}\right)^{\frac{1}{s_{k}}}\right)
-2\pi
p_{0}dI\left(\dfrac{1}{\varepsilon_{n}}\left(\log\dfrac{1}{\varepsilon_{n}}\right)^{-\frac{1}{s_{k}}}\right)\\
&+ {9\over 8}\pi p_{0}\sum_{k=1}^{N}\sum_{i\neq j}\log \dfrac{T_{\varepsilon_{n}}}{\mid y_{i}-y_{j}\mid}+O(1)
\end{split}
\end{equation}
which is $\left(\ref{lb}\right)$.
\begin{remark}\label{remDPF}
In Proposition \ref{proplb} we have proved \eqref{lb} for a starshaped domain. An argument of del Pino and Felmer in \cite{PF} can now be used to show that \eqref{HS3.2} holds without the assumption on the starshapeness of $G$. Hence \eqref{lb} is still true for general domain and we can conclude again by acting as in \cite{HS1}.
\end{remark}
\subsection{Conclusions}\label{secend}
By collecting together Proposition $\ref{upb}$ and Proposition $\ref{proplb}$, and taking into account Remark \ref{remDPF}, we obtain Proposition \ref{as} which is \eqref{conv2} of Theorem \ref{teo1}.

Thanks to estimate \eqref{HS3.2}, we can now follow the construction of bad discs as in \cite{BBH1} and prove convergence \eqref{conv1} of Theorem \ref{teo1}. Since the arguments are identical to those of \cite{BBH1} we omit the details. Now Theorem \ref{teo1} is completely proved.\\
Finally as a consequence of  $(\ref{up10})$ and $(\ref{lbk})$, we get the following  estimate of the distance between the centers of bad discs. 

\begin{corollary}\label{propAS.3.1} 
For every $i \neq j$ in $I_{k}$ $(1\leq k\leq N_{2})$ with $\left|I_{k}\right|=d_{k}>1$, we have
\begin{equation}\label{AS3.13*}
   exp\left( - C_{8} I \left(\left(\log\dfrac{1}{\varepsilon_{n}}\right)^{\frac{1}{s_{k}}}\right)\right)  \mid\log\varepsilon_{n}\mid^{-\frac{1}{s_{k}}} \leq  \left|y_{i}-y_{j}\right | \leq C_{9}\mid\log\varepsilon_{n}\mid^{-\frac{1}{s_{k}}}
\end{equation}
where $C_8$ and $C_9$ are two constants independent of $\varepsilon$.
\end{corollary}
\textbf{Proof} 
By lower bound $(\ref{lbk})$ we have
\begin{equation}\label{BMR5***} 
\begin{split}
   \int_{\Omega}p\left|\nabla u\right|^{2}\geq &
 2\pi p_{0}d_{k}\log\dfrac{1}{\varepsilon_{n}} +   2\pi p_{0}\dfrac{d_{k}^{2}-d_{k}}{s_{k}}\log\log\dfrac{1}{\varepsilon_{n}}-2\pi
    p_{0}\Sigma_1^N d_{k}^{2} I \left(\left(\log\dfrac{1}{\varepsilon_{n}}\right)^{\frac{1}{s_{k}}}\right)
    \\&
     -2\pi
p_{0}d I\left(\dfrac{1}{\varepsilon_{n}}\left(\log\frac{1}{\varepsilon_{n}}\right)^{-\frac{1}{s_{k}}}\right)+{9\over 8}\pi p_{0}\sum_{i\neq j}\log \dfrac{T_{\varepsilon_{n}}}{\mid y_{i}-y_{j}\mid}+O(1).
\end{split}
\end{equation}

The upper bound $(\ref{up10})$ and \eqref{BMR5***}, imply
\begin{equation*}
\sum_{i\neq j}\log \left(\dfrac{\mid\log\varepsilon_{n}\mid^{-\frac{1}{s_{k}}}}{\mid y_{i}-y_{j}\mid}\right)\leq 
 C_{8}I \left(\left(\log\dfrac{1}{\varepsilon_{n}}\right)^{\frac{1}{s_{k}}}\right)
\end{equation*}
which, by using \eqref{AS3.13},
is the claimed result.

\end{document}